\documentclass[reqno,11pt]{amsart}
\usepackage{amsfonts,amsmath,amssymb,amsthm,mathtools}
\usepackage{url}
\usepackage[colorlinks=true,citecolor=blue]{hyperref}
\usepackage[numbers,sort&compress]{natbib}
\usepackage[toc,page]{appendix}
\usepackage{latexsym}
\usepackage{hyperref}
\usepackage{enumerate}
\RequirePackage{natbib}
\usepackage{natbib}
\usepackage{import}
\usepackage{xcolor}
\usepackage[utf8]{inputenc}
\setcitestyle{square}
\usepackage{url}
\usepackage{algorithm}
\usepackage[noend]{algpseudocode}
\usepackage{multirow} 
\usepackage{graphicx}
\usepackage{subfig}
\usepackage{epstopdf}
\usepackage{pgfplots}
\usepackage{pgfplotstable}
\usepackage[pdf]{pstricks}
\usepackage{tikz}
\usetikzlibrary{arrows,positioning,chains,fit,shapes,calc,snakes}
\usepackage{lipsum}

\newtheorem{theorem}{Theorem}[section]

\newtheorem{proposition}[theorem]{Proposition}

\newtheorem{conjecture}[theorem]{Conjecture}

\def\RR{{\mathbb R}}  
\def\PP{{\mathbb P}}
\def\NN{{\mathbb N}}      
\def\Fix{\mathop{\rm Fix}}      
 
\DeclareMathOperator{\diam}{diam}
\def\Pcal{{\mathcal P}}
\def\Fcal{{\mathcal F}}
\def\Ropt{{\hat{R}}}

\def\D{\mathcal{I}}

\title[Optimal bounds  for nonexpansive fixed point iterations]{Optimal error bounds for nonexpansive fixed-point iterations in normed spaces}

\author[J.P. Contreras]{Juan Pablo Contreras}
\address[J.P.C.]{Universidad Adolfo Ib\'a\~nez, Facultad de Ingenier\'ia y Ciencias, Diagonal Las Torres 2640, Santiago, Chile}
\email{ \href{mailto:juan.contrerasff@gmail.com}{\nolinkurl{juan.contrerasff@gmail.com}}}
\thanks{} 

\author[R. Cominetti]{Roberto Cominetti}
\thanks{}
\address[R.C.]{Universidad Adolfo Ib\'a\~nez, Facultad de Ingenier\'ia y Ciencias, Diagonal Las Torres 2640, Santiago, Chile} 
\email{ \href{mailto:roberto.cominetti@uai.cl}{\nolinkurl{roberto.cominetti@uai.cl}}}

\begin{document}
\bibliographystyle{plainnat}

 \maketitle

\begin{abstract}
This paper investigates optimal error bounds and convergence rates for general Mann iterations for computing fixed-points of non-expansive maps. We look for iterations that achieve the smallest fixed-point residual after $n$ steps, by  minimizing a worst-case bound $\|x^n-Tx^n\|\le R_n$ derived from a nested family of optimal transport problems. We prove that this bound is tight so that minimizing $R_n$ yields optimal iterations. Inspired from numerical results we identify iterations that attain the rate $R_n=O(1/n)$, 
which we also show to be the best possible.
In particular, we prove that the classical Halpern iteration achieves this optimal rate
for several alternative stepsizes, and we determine analytically the optimal stepsizes that 
attain the smallest worst-case residuals at every step $n$, with a tight bound $R_n\approx\frac{4}{n+4}$. 
We also determine the optimal Halpern stepsizes for
affine non-expansive maps, for which we get exactly $R_n=\frac{1}{n+1}$. Finally, we show that the best rate for the classical Krasnosel'ski\u{\i}-Mann iteration is 
$\Omega(1/\sqrt{n})$, and present numerical evidence suggesting that even extended variants cannot reach a faster rate.
\end{abstract}

\begin{small}
  \noindent{\bf Keywords:} Nonexpansive maps $\cdot$ fixed-point iterations $\cdot$ error bounds $\cdot$ convergence rates \\[2ex]
  \noindent{\bf  Mathematics Classification Subjects 2020:} 47J25 $\cdot$ 47J26 $\cdot$ 65J15 $\cdot$ 65K15 \\[2ex]
  \noindent Forthcoming in Mathematical Programming (2022)
  \end{small}


\section{Introduction}

Many computational problems can be cast as finding a fixed point of a  map
$T:C\mapsto C$ where $C\subseteq X$ is a bounded convex domain on a normed space $(X,\|\cdot\|)$
and $T$ is non-expansive, that is
\begin{equation}
\|Tx-Ty\|\le \|x-y\|, \qquad \mbox{for all } x,y\in C.
\label{non-expansivedef}
\end{equation}
When $T:X\to X$ is  defined on the full space and has a fixed point $x^*=Tx^*$, 
one can take $C=B(x^*,r)$ as any ball  centered at $x^*$ with radius $r\geq 0$ and $\diam(C)=2 r$. 
In general,  by
rescaling the norm we may assume without  loss of generality 
that $\mathop{\rm diam}(C)=1$, which we do from now on.

Fixed-point iterations arise in different settings, including, among others, decomposition methods in convex optimization, regression in statistical estimation, computation of invariant measures of Markov chains, solution of monotone inclusions,
asymptotics of dissipative dynamical systems, and more. Consequently, there is a variety of iterative
methods for solving such problems. 
When $X$ is a Banach space and $T$ is a strict contraction, the method of choice
is the classical Banach-Picard iteration $x^{n}=Tx^{n-1}$ which converges at an $R$-linear rate to the unique fixed point 
of $T$. If $T$ is just non-expansive 
this is no longer true: the map may not only fail to have fixed points but even if $\Fix(T)$ is nonempty the
Banach-Picard iterates might not converge. 
A classical strategy to overcome this issue is to consider an {\it averaged mapping} ({\it cf.} \cite{bailion1978asymptotic}), 
$T_\alpha = (1-\alpha)I+\alpha T$ with  $\alpha\in (0,1)$, which has the same fixed points as $T$.
Krasnosel'ski\u{\i} \cite{krasnosel1955two} considered the case $\alpha\equiv 1/2$ and established 
that the Banach-Picard iteration applied to $T_\alpha$ produces a sequence $\{x^n\}_{n\in\mathbb{N}}$ 
that converges in norm to a fixed point provided that the space $X$ is uniformly convex and $T(C)$ is relatively 
compact. Under weaker conditions, one can still obtain weak convergence to a fixed point of $T$ (see Reich \cite{reich1979} and Borwein et. al. \cite{borwein1992krasnoselski}). By considering different coefficients 
$\alpha_n$ at each step we obtain the Krasnosel'ski\u{\i}-Mann iteration $x^{n} = (1-\alpha_n)x^{n-1}+\alpha_nTx^{n-1}$ for which similar convergence results hold. 
Also, replacing $x^{n-1}$ in this average with a constant vector $y^0$ we get the classical Halpern iteration $x^{n} = (1-\alpha_n)y^{0}+\alpha_nTx^{n-1}$ ({\em cf. } \cite{Halpern1967fixed}). Mann \cite{mann1953mean} proposed a more general iterative scheme that involves a weighted average of all the previous iterates,
which was also used by Reich \cite{reich1975fixed} in order to study the existence of fixed points. For a comprehensive survey of applications, methods,  and convergence results, we refer to Berinde \cite{berinde},  and  Bauschke \& Combettes \cite{bauschke2011convex}.

Most of the existing methods turn out to be special cases of the general Mann iteration
recalled in  \eqref{generalmann0} below. Our main goal  
is to investigate this general framework in order to design fast iterations that attain the smallest possible norm for 
the fixed-point residual $\|x^n-Tx^n\|$ after $n$ steps. To this end, we exploit some universal bounds
$\|x^n-Tx^n\|\leq R_n$ where $R_n$ is obtained through a nested sequence of optimal transport problems. 
Along the way, we establish some lower bounds on what can be achieved.

\begin{figure}[t]
\centering
\includegraphics[scale=1]{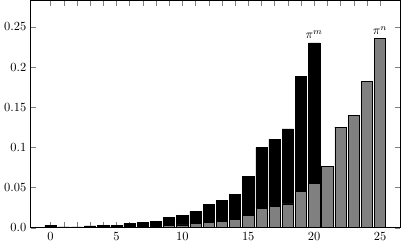}
\caption{\label{pis}The distributions $\pi^m$ and $\pi^n$ ($m=20, n=25$).}
\end{figure}

Given a pair of initial points $x^0,y^0\in C$, and adopting the convention $Tx^{-1}=y^0$, Mann's iterates are  defined recursively as
\begin{equation}\label{generalmann0}
x^n = \sum_{i=0}^n \pi^n_i\,Tx^{i-1} \qquad (\forall\,n\ge 1)
\end{equation}
where $\pi^n = (\pi^n_i)_{i=0}^n$ is a given sequence of averaging coefficients with $\pi^n_i\geq 0$ and $\sum_{i=0}^n\pi_i^n = 1$.
Such $\pi^n$ can be seen as a probability distribution on $\NN$ supported on $\{0,1,\ldots,n\}$ (see Figure \ref{pis}).
A common choice is to take $y^0=x^0$, however  it is sometimes convenient to take $y^0\neq x^0$
as in Halpern's iteration.

Mann's iteration \eqref{generalmann0} is very general and includes, among others, the classical Krasnosel'ski\u{\i}, 
Halpern, and Ishikawa iterations. More recently, motivated by the 
accelerated gradient methods in smooth convex optimization that incorporate inertial terms ({\em cf.} Nesterov \cite{nesterov2018}), there 
has been a renewed interest in adapting these  techniques in order to accelerate fixed-point iterations. 
These classical iterations and some modifications including extra terms in the averaging will be analyzed in Section \ref{specialMann}.

\subsection{Previous results on convergence rates}
A crucial step in establishing the convergence of the iterates is to bound the 
fixed point residual $\|x^n-Tx^n\|$ and to prove that it converges to 0, a property named as {\it asymptotic regularity}   
after Browder \& Petryshyn \cite{browder1966solution}. For descent methods in smooth convex optimization, this corresponds to the 
strong convergence of the gradient $\|\nabla f(x^n)\|\to 0$.

For the Krasnosel'ski\u{\i}-Mann iteration $x^{n+1}=(1\!-\!\alpha_n)x^n+\alpha_n Tx^n$, convergence rates were first obtained by Baillon \& Bruck \cite{baillon1992optimal,baillon1996rate} for 
constant stepsizes  $\alpha_n\equiv\alpha$, and later  extended in Cominetti {\em et al.} \cite{cominetti2014rate} establishing the estimate
$$\|x^n-Tx^n\|\leq\frac{\diam(C)}{\sqrt{\pi \sum_{i=1}^n\alpha_i(1\!-\!\alpha_i)}}.$$
For constant stepsizes $\alpha_n\equiv\alpha$ this yields the rate $\|x^n-Tx^n\|\sim O(1/\sqrt{n})$.
Bravo \& Cominetti \cite{bravo2018sharp}  later proved  that the constant $1/\sqrt{\pi}$ in this bound is tight
and cannot be improved in general normed spaces. 

A natural question, which was only recently settled, is whether a faster rate $O(1/n)$ could be achieved 
with suitably chosen $\pi^n$'s. A positive answer came from another special case of 
\eqref{generalmann0},
namely, the classical Halpern iteration $x^{n+1}=(1\!-\!\beta_n)y^0+\beta_n Tx^n$.
Halpern \cite{Halpern1967fixed} considered the special case  with $C$ the unit ball of a Hilbert space and $x^0=y^0=0$,
establishing necessary and sufficient conditions on $\beta_n$ to ensure the strong convergence towards
a fixed point. Wittmann \cite{wittmann1992approximation} supplemented this result by proving strong convergence 
to the fixed point closest to $y^0$,{ provided that $\beta_n\rightarrow 1$, $\sum_n (1-\beta_n)=+\infty$, and 
$\sum_n\vert\beta_{n+1}-\beta_n\vert < +\infty$. This was the first convergence result that covered the classical 
coefficients $\beta_n = \frac{n}{n+1}$.  In uniformly smooth Banach spaces, Reich \cite{reich1980} proved strong convergence for the particular choice $\beta_n = 1-\frac{1}{(n+2)^a}$ with $0<a<1$, and extended Wittmann's theorem in \cite{reich1994}. Using proof-mining techniques, the works of Leustean \cite{leustean2007rates} and 
Kohlenbach \cite{kohlenbach2011quantitative} were able to extract convergence rates for Halpern's iteration from the proofs in Browder \cite{browder1967convergence} and Wittmann \cite{wittmann1992approximation}. Uniform rates 
of metastability were proved by Kohlenbach \cite{kohlenbach2010logical,kohlenbach2011quantitative}, while the paper
K\"ornlein  \cite{kornlein2015quantitative} established a rate of metastability extracted from a proof in 
 Xu \cite{xu2002iterative} in the setting of uniformly smooth spaces.

To the best of our knowledge, the first proof of a rate $O(1/n)$ for Halpern's iteration was established
in Sabach \& Shtern \cite{sabach2017first}. Using the stepsizes $\beta_n=\frac{n}{n+2}$ they found the explicit
bound $\|x^n-Tx^n\|\leq\frac{4}{n+1}$ (see \cite[Lemma 5]{sabach2017first}). Independently, Lieder \cite{lieder2021convergence}\footnote{First version appears in 2017 in
\url{optimization-online.org/DB_FILE/2017/11/6336.pdf}} studied 
Halpern's iteration in Hilbert spaces by using techniques of Performance Estimation Problems (PEP) (Drori \& Teboulle \cite{drori2014performance}), and established that $\beta_n=\frac{n}{n+1}$ improves this estimate to  $\|x^n-Tx^n\|\leq\frac{1}{n+1}$,
providing also a simple example where this bound is attained. Hence, the classical Halpern iteration  not only  
achieves a sharp bound but seems to be optimal among all Mann-type 
iterations in Hilbert spaces. A formal proof of the latter is still an open question. In a similar direction, Kim
 \cite{kim2021accelerated}\footnote{First version appears in 2019 in {\it arXiv
preprint arXiv:1905.05149}} proposed a general Mann-type algorithm to find a zero of  
a co-coercive operator. By  numerically computing the optimal coefficients using PEP techniques, he discovered an inertial iteration that
achieves the same bound as Lieder's. As a matter of fact, we will show that when
translated into the setting of fixed points for non-expansive maps, the inertial iteration in 
Kim \cite{kim2021accelerated} coincides with a classical Halpern iteration and the result is equivalent to
Lieder's \cite{lieder2021convergence}. Finally, still in Hilbert spaces, Diakonikolas \cite{diakonikolas2020Halpern} provides a simpler 
potential-based proof of the rate $O(1/n)$ which applies to more general stepsizes and motivates parameter-free algorithms 
for monotone inclusions, variational inequalities, convex-concave min-max optimization, and related problems.

\subsection{Our contribution}
In this paper we investigate the convergence rate of general Mann iterations for non-expansive maps in normed spaces. 
Building upon the optimal transport approach from Cominetti {\em et al.} \cite{cominetti2014rate}, we formulate an optimization problem that aims to find 
the sequence of averaging parameters $\pi^n$ that minimize a worst-case  bound $R_n=R_n(\pi)$ for the fixed point residual $\|x^n-Tx^n\|\leq R_n$. 
We prove that this bound is tight so that minimizing $R_n(\pi)$ yields the best possible iteration. We show that the optimal rate of convergence for 
general Mann iterations is $O(1/n)$, presenting a simple linear map $T$ for which the residual in every Mann iteration satisfies $\|x^n-Tx^n\|\geq \frac{1}{n+1}$.

By incorporating additional structural constraints in the parameters $\pi^n$, we investigate some special cases of the  Mann iteration such as Krasnosel'ski\u{\i}-Mann, Halpern, Ishikawa, as well as some variants  that include extra terms in the averaging. Through numerical computations, we identify several simple iterations that attain the optimal rate $O(1/n)$. 
For Halpern's iteration we 
 prove  that the optimal stepsizes are given by the recursion $\beta_{n+1} = (1+\beta_{n}^2)/2$ with $\beta_0 = 0$. These stepsizes yield a tight bound 
$\|x_n-Tx_n\|\le R_n$ that satisfies $R_n\leq \frac{4}{n+4}$, slightly improving upon Sabach \& Shtern \cite{sabach2017first}. 

We also revisit Halpern's iteration in Hilbert spaces, and highlight the connections between the recent papers by
Lieder \cite{lieder2021convergence} and Kim \cite{kim2021accelerated} which establish the tight bound $\frac{1}{n+1}$ for general non-expansive maps. Remarkably, we show that for {\em affine} non-expansive maps in Banach spaces the bound $R_n=\frac{1}{n+1}$ is also tight, and moreover the optimal coefficients  in both settings coincide. 

Finally, we investigate the optimal stepsizes for the Krasnosel'ski\u{\i}-Mann iteration, as well as some variants that have 
received attention recently. We show that the best possible rate in any Krasnosel'ski\u{\i}-Mann iteration is $O(1/\sqrt{n})$,
a rate that is known to be achieved with constant stepsizes ({\em cf.} Baillon \& Bruck \cite{baillon1996rate}).
Our numerical results suggest that the variants of Krasnosel'ski\u{\i}-Mann including extra terms do not improve this rate, even when using optimal averaging coefficients.

\subsection{Structure of the paper} 
Section \ref{generalmanniterationsection} recalls the general setting of Mann's iterations
and presents the optimal transport bounds for the fixed point residuals. After summarizing previously known properties 
and establishing the tightness of the worst-case bound, we perform a numerical optimization of the general Mann
iterations as well as several special subclasses. Motivated by these numerical results, Section \ref{LB} establishes
explicit lower bounds for the best rates that can be achieved by Mann as well as Krasnosel'ski\u{\i}-Mann iterations, after which 
Section \ref{optHalpern} presents a detailed analysis of the optimal Halpern iteration and several variants that attain the 
optimal convergence rate $O(1/n)$.

\section{General Mann iterations}
\label{generalmanniterationsection}
This section describes a framework for finding the optimal averaging parameters for 
Mann's fixed points iterations \cite{mann1953mean}. 
We recall that these iterations are defined by a triangular array 
$$\pi  = \left(\begin{array}{cccccc}
\pi^0_0 &  & & & &\\
\pi^1_0 & \pi^1_1 &  & & &\\
\pi^2_0 & \pi^2_1 & \pi^2_2 &  & &\\
\vdots & \vdots & \vdots & \ddots &  &\\
\pi^n_0 & \pi^n_1 & \pi^n_2 & \cdots & \pi^n_n&\\
\vdots & \vdots & \vdots & & \vdots &\ddots \\
\end{array}\right)$$
where the $n$-th row $\pi^n$ satisfies $\pi^n_i\ge 0$, $\sum_{i=0}^n \pi^n_i = 1$, and $\pi_i^n=0$ for $i>n$. Note that  $\pi^n$ can be interpreted as a probability distribution on  the set of non-negative integers  $\NN = \{0,1,2,...\}$, with support included in $\{0,1,\ldots,n\}$. 

Given initial points $x^0,y^0\in C$ and adopting the convention $Tx^{-1}=y^0$,
the Mann iterates are given by the following sequential averaging scheme
\begin{equation}\label{generalmann}
 x^n = \sum_{i=0}^n \pi^n_i\, Tx^{i-1}\qquad (\forall\,n\ge 1).
\end{equation}
This iteration is very general and includes the classical Krasnosel'ski\u{\i}-Mann, 
Halpern, and Ishikawa iterations (see Section \ref{specialMann}). We are interested in determining the matrix of coefficients $\pi$ that minimize the norm of the residuals $\|x^n-Tx^n\|$ after $n$ iterations. More precisely,
we seek to minimize the worst-case among all possible runs of \eqref{generalmann}, that is
\begin{equation}\label{worstResidual}
\Psi_n(\pi)=\sup_{T,x^0\!,y^0}\|x^n-Tx^n\|
\end{equation}
where the supremum is taken over all non-expansive maps $T:C\to C$ with $\diam(C)=1$, 
and all possible initial points $x^0,y^0\in C$.

\subsection{Sharp error bounds via optimal transport metrics}\label{optbounds}
The expression  \eqref{worstResidual} for the worst-case error bound  is not easy to handle. 
Below we present an alternative formula $\Psi_n(\pi)=R_n(\pi)$ which is more manageable in view of optimizing 
over $\pi$. 
In order to estimate the fixed point residuals $\|x^n-Tx^n\|$ we follow the approach in Bravo \& Cominetti \cite{bravo2018sharp} 
and further developed in Bravo {\em et al.} \cite{bcc2020}. This
is based on estimates of the distance between iterates $\|x^m-x^n\|\leq d_{m,n}$, where
 the bounds $d_{m,n}$ are defined by a nested sequence of optimal transport problems.

Starting with  $d_{-1,-1}=0$ and $d_{-1,j}=d_{j,-1}= {\rm diam}(C)=1$ for all $j\in\NN$, we consider 
the double-indexed family of reals $d_{m,n}$ for $m,n\in\NN$ defined by 
\begin{equation*}\label{dmn}
d_{m,n}=\min_{z\in\Fcal(\pi^m,\pi^n)} \sum_{i=0}^m\sum_{j=0}^n z_{i,j}\,d_{i-1,j-1}\leqno (\Pcal_{m,n})
\end{equation*}
where $\Fcal(\pi^m,\pi^n)$ denotes the set of transport plans taking $\pi^m$ into $\pi^n$, that is,
the set of all $z=(z_{i,j})_{i=0,\ldots,m; j=0,\ldots n}$ such that $z_{i,j}\geq 0$ and
$$\begin{array}{ll}
\sum_{j=0}^n z_{i,j}=\pi_i^m&\mbox{for all } i=0,\ldots,m;\\[1.5ex]
\sum_{i=0}^m z_{i,j}=\pi_j^n&\mbox{for all } j=0,\ldots,n.
\end{array}
$$

Each transport plan $z$ from $\pi^m$ to $\pi^n$ yields the estimate
\begin{eqnarray}\nonumber
\label{xmn}
\|x^m-x^n\|&=&\mbox{$\|\sum_{i=0}^m\sum_{j=0}^nz_{i,j}(Tx^{i-1}-Tx^{j-1})\|$}\\ \nonumber
&\leq&\mbox{$\sum_{i=0}^m\sum_{j=0}^nz_{i,j}\|Tx^{i-1}-Tx^{j-1}\|$}\\ \nonumber
&\leq&\mbox{$\sum_{i=0}^m\sum_{j=0}^nz_{i,j}d_{i-1,j-1}$}
\end{eqnarray}
where the last inequality uses the non-expansivity of $T$  and assumes that we already have 
$\|x^{i-1}-x^{j-1}\|\leq d_{i-1,j-1}$ for previous iterates (for $i=0$ or $j=0$ use the 
bound $\|y^{0}-Tx^{k}\|\leq \mathop{\rm diam}(C)=1$).
Minimizing over $z$ we get 
\begin{equation}\label{bound_dmn}
\|x^m-x^n\|\leq d_{m,n}
\end{equation}
from which it follows inductively that this inequality holds for all $m,n\in\NN$. 
Using the triangle inequality and non-expansivity, we can  then
estimate the fixed-point residuals as
\begin{eqnarray}\nonumber
\|x^n-Tx^n\|&=&\mbox{$\|\sum_{i=0}^n\pi^n_i(Tx^{i-1}-Tx^n)\|$}\\
&\leq& \mbox{$\sum_{i=0}^n\pi_i^nd_{i-1,n}\triangleq R_n$}.
\label{bound_Rn}
\end{eqnarray}

We emphasize that the optimal transport bounds $d_{m,n}=d_{m,n}(\pi)$ and $R_n=R_n(\pi)$ are universal in the sense that they  only depend on the sequence 
$\pi$ and not on the particular map $T$ or the initial points $x^0$ and $y^0$, so that 
\begin{equation}\label{worstResidual2}
\Psi_n(\pi)\leq R_n(\pi).
\end{equation}
In fact, for $m\leq n$ both 
$d_{m,n}$ and $R_n$ only depend on $\pi^0,\pi^1,\ldots,\pi^n$. Note also that, by symmetry it suffices to compute 
$d_{m,n}$ for $m\leq n$.

It turns out that these bounds are tight.
The proof is based on ideas that evolved from the original paper by Baillon \& Bruck \cite{baillon1992optimal},
and  further developed in Bravo \& Cominetti \cite{bravo2018sharp} and Bravo {\em et al.} \cite{bcc2020}. 
These latter references construct non-expansive 
maps $T$ and sequences $(x^n)_{n\in\NN}$ that attain all the bounds 
\eqref{bound_dmn} and \eqref{bound_Rn}  with equality, first for Krasnosel'ski\u{\i}-Mann iterations \cite{bravo2018sharp}, and later
in \cite{bcc2020}
for every Mann iteration satisfying the additional monotonicity condition
 \begin{equation}\label{eq_monotone}
 (\forall n\geq 1)\quad \mbox{$\pi_n^n>0$ and $0\leq \pi_i^{n}\leq\pi_i^{n-1}$ for $i<n$}.
 \end{equation}
Here we take one step forward by showing that the bounds \eqref{bound_dmn} and \eqref{bound_Rn} 
are always the best possible and are attained with equality for a suitably chosen map, without any extra condition 
such as \eqref{eq_monotone}.

\begin{theorem}\label{TightMann}
Let $\pi=(\pi^n)_{n\in\NN}$ be a sequence such that $\pi^n_i\ge 0$, $\sum_{i=0}^n \pi^n_i = 1$, and $\pi_i^n=0$ for $i>n$. Then there exists a non-expansive map $T$ and a corresponding Mann sequence $(x^n)_{n\in\NN}$ that attains all the bounds \eqref{bound_dmn} and \eqref{bound_Rn} with equality. In particular $\Psi_n(\pi)=R_n(\pi)$ for all $n\in\NN$.
\end{theorem}
\begin{proof} See Appendix \ref{SecTightMann}. \end{proof}

The optimal transport bounds will be used to design Mann iterations that 
minimize the worst-case resisual bounds $\Psi_n(\pi)=R_n(\pi)$.
Before proceeding we summarize the main results from Bravo {\em et al.} \cite{bcc2020} that will be used later. 
To begin with, it was shown that $(m,n) \mapsto d_{m,n}$ defines a metric on the set $N\triangleq \{-1,0,1,2,\ldots\}$ 
with $d_{m,n}\in[0,1]$ for all $m,n\in N$. This implies the alternative dual characterization
\begin{equation*}
\begin{array}{cl}
d_{m,n}=&\displaystyle{\max_{u}}~~\mbox{$\sum_{j=0}^{n}$}\,(\pi^n_j\! -\! \pi^m_j)\,u_j\\[1ex]
&\mbox{s.t.}~\vert u_i-u_j\vert\leq d_{i-1,j-1}\mbox{ $\forall i,j=0,\ldots,n$}
\end{array} \leqno (\mathcal{D}_{m,n})
\end{equation*}
so that each pair of primal-dual optimal solutions $z^{mn}$ and $u^{mn}$ 
satisfy the complementary slackness 
 \begin{equation}\label{cs}
z_{i,j}^{mn}(u_j^{mn}\!-u_i^{mn})=z^{mn}_{i,j}d_{i-1,j-1}\mbox{ for all $i,j=0,\ldots,n$}.
\end{equation}
Using these facts the bound \eqref{bound_dmn} was proved to be tight, showing also that
each $(\Pcal_{m,n})$ has a {\em simple} optimal transport $z^{m,n}$ with
$z^{m,n}_{i,i}=\min\{\pi^m_i,\pi^n_i\}$ for all $0\leq i\leq\min\{m,n\}$. 
Moreover, under the monotonicity condition \eqref{eq_monotone}  it was proved
that \eqref{bound_Rn}  is also tight, and that the $d_{m,n}$'s satisfy 
the convex quadrangle inequality (see Figure \ref{fig:4points})
 \begin{equation}\label{4point}
d_{i,l}+d_{j,k}\leq d_{i,k}+d_{j,l}\qquad   \mbox{for all $i<j<k<l$.}
\end{equation}
The latter yields a greedy algorithm to compute optimal transports
that are {\em nested} in the sense that the flows do not intersect. 

\begin{figure}[t]
\centering
\includegraphics[scale=1]{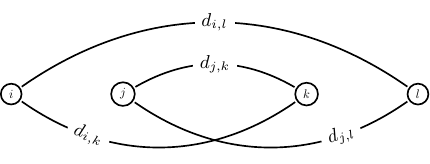}
\caption{\label{fig:4points} The convex quadrangle inequality.\label{4-point}}
\end{figure}

\subsection{Fixed-horizon optimization of Mann iteration}
The  optimal transport bounds can be exploited to 
determine the averaging coefficients $\pi$ that yield the smallest possible
worst-case residual $\Psi_n(\pi)=R_n(\pi)$  for any fixed horizon $n$.
This can be stated as the non-convex optimization problem
\begin{equation*}
\begin{array}{rll}
\mbox{\sc (fh$_n$)}\hspace{0.5cm} \displaystyle{\min_{d,z,\pi}} & \sum_{k=0}^n\,\pi^{n}_k \,d_{k-1,n} &\\
\mbox{s.t.} 
& d_{-1,k} = 1, &\forall k=0,...,n\\
& d_{k,k} = 0, &\forall k=-1,...,n\\
& d_{k,m} = {\sum_{i=0}^{k}\sum_{j=0}^{m}} \,z^{k,m}_{i,j}\,d_{i-1,j-1}, &\forall 0\leq k<m\leq n\\
& d_{m,k} = d_{k,m}, &\forall 0\leq k<m\leq n\\
& \pi^k\in \Delta^k, \quad z^{k,m} \in \Fcal(\pi^k,\pi^m), &\forall 0\leq k<m\leq n.
\end{array}
\label{proboptmann}
\end{equation*}
where $\Delta^k=\{(x_0,\ldots,x_k):x_i\geq 0, \sum_{i=0}^k x_i=1\}$ denotes the unit simplex in dimension $k+1$.

The objective function is precisely the bound $R_n(\pi)$ at the $n$-th iteration. 
The first constraint sets the boundary conditions for the distances $d_{-1,k}$, while
the second to fourth constraints correspond to the successive use of the optimal transport
problems $(\Pcal_{k,m})$ that determine the bounds for the distance between iterates.  
In order to see that this is a valid formulation and that these constraints induce the 
variables $z^{k,m}$ to be optimal transports, we observe that the objective function seeks to make $d_{k-1,n}$ as small as
possible. This fact, combined with the recursive and monotonic dependence of $d_{k-1,n}$
on previous $d_{i,j}$'s, implies that all these variables should be made small. This
automatically pushes the variables $z^{k,m}$ to be  chosen 
as optimal transports that minimize the right-hand side in the third equality constraint.

Problem {\sc (fh$_n$)} turned out to be quite hard to solve. Experimental computations 
with state-of-the-art solvers show poor performance already for $n\ge 7$. The difficulties arise from 
the large number of variables and the non-linearities both in the objective function as well as in the 
constraints which involve multiple nested products of the variables $\pi, d$ and $z$. 
A further drawback of this global optimization approach is the fact that the optimal solutions of \
{\sc (fh$_n$)} depend on the pre-fixed horizon $n$, and the full optimal sequence 
$\pi^0,\pi^1,\ldots,\pi^n$ is modified when $n$ changes. This is inconvenient if one does not 
know a priori the number of iterations to be performed. This leads to consider a 
simplification in which the $\pi^n$'s are optimized sequentially by fixing the previous solutions,
which also makes the subproblems computationally tractable.

\subsection{Sequential optimization of Mann iteration}
In the sequential optimization approach we seek to determine the averaging coefficients $\pi$ 
without relying on a pre-established number of iterations.  Specifically, we compute $\pi^n$ progressively
by using the optimization problem {\sc (fh$_n$)} but fixing the solutions $\pi^0,...,\pi^{n-1}$ and distances $d_{i,j}$ 
computed in the previous stages $1,...,n\!-\!1$.  By considering the 
optimal transport bounds $d_{k,n}=d_{k,n}(\pi^n)$ as a function of $\pi^n$, with $\pi^0,\ldots,\pi^{n-1}$ fixed,
the $n$-th stage problem becomes 
\begin{equation*}
\min_{\pi^n \in \Delta^{\!n}} R_n(\pi^n)\triangleq \sum_{k=0}^n\,\pi^n_k \, d_{k-1,n}(\pi^n)
\end{equation*} 
or more explicitly 
\begin{equation*}
\begin{array}{rll}
\mbox{\sc (s$_n$)}\hspace{0.5cm}\displaystyle{\min_{\pi^n,z}} & \pi^n_0+\displaystyle{\sum_{k=1}^n} \sum_{i=0}^{k-1}\sum_{j=0}^n\pi^{n}_{k}z^{k-1,n}_{i,j}d_{i-1,j-1} &\\
\mbox{s.t.} & & \\
 & \pi^n\in \Delta^{\!n}, \quad z^{k,n} \in \Fcal(\pi^k,\pi^n),& \forall k=0,...,n-1.
\end{array}
\label{mannasymptoticaux}
\end{equation*} 

These problems are still non-convex but the number of variables is significantly reduced. In particular 
the distances $d_{i-1,j-1}$ in the objective function  are fixed and are no longer variables. As a consequence 
$(S_n)$ becomes a  linearly constrained (non-convex) quadratic programming problem which can 
be solved with general non-linear solvers such as Knitro or Baron, or extended linear programming solvers such as CPLEX and Gurobi.

\subsection{Monotone sequential optimization of Mann's iteration}
\label{simpletransport}

A further simplification of the sequential approach is achieved by restricting the $\pi^n$'s to satisfy the 
monotonicity condition \eqref{eq_monotone}. As mentioned earlier, in this case 
we have the convex quadrangle inequality and there is a simple greedy algorithm to compute the optimal 
transports. These optimal transports have a particularly simple structure when we have in addition 
\begin{equation}\label{eqm2}
 \pi_{m}^m\ge \mbox{$\sum_{i=m}^{n-1} \pi^n_{i},\qquad\forall m<n$}.
 \end{equation}
 Note that this holds automatically if $\pi^m_m$ and $\pi^n_n$ are larger than $\frac{1}{2}$.
The following result presents the simpler explicit expression for the $n$-th stage residual, which avoids the
use of the $z$ variables. This can be 
derived from results in \cite{bcc2020}. Here we present an alternative proof using linear programming 
duality. 
\begin{proposition}
\label{propsubproblem}
Assume the quadrangle inequality \eqref{4point} holds up to $n-1$.
Let $m<n$ and suppose further \eqref{eqm2} and $\pi_i^n\leq\pi_i^m$ for $i=0,\ldots,m$. 
Then, $d_{m,n}(\pi^n)=D_{m,n}(\pi^n)$ where 
\begin{equation}
D_{m,n}(\pi^n)\triangleq \sum_{i=0}^m (\pi^m_i\!-\pi_i^n)d_{i-1,n-1}+\!\!\!\sum_{j=m+1}^n\!\!\!\pi_j^n(d_{m-1,j-1}-d_{m-1,n-1}).
\label{optvaluedmn}
\end{equation}

\end{proposition}

\begin{proof} A straightforward computation shows that $D_{m,n}(\pi^n)$ is 
exactly the cost of 
 the following feasible transport plan for $d_{m,n}$  
\begin{equation}
\begin{array}{ll}
z_{i,i} = \pi^n_i & \forall i=0,...,m \\
z_{m,j} = \pi^n_j & \forall j=m+1,...,n-1 \\
z_{i,n} = \pi^m_i-\pi_i^n & \forall i=0,...,m-1\\
z_{m,n} = \pi^m_m-\sum_{j=m}^{n-1}\pi_j^n .& 
\end{array}
\label{zstructure}
\end{equation}
Let us consider the dual of the linear program $(\Pcal_{m,n})$ given by
$$
\begin{array}{cl}
\max_{u,v}& \sum_{j=0}^n \pi^n_j u_j-\sum_{i=0}^m \pi^m_i v_i\\
\mbox{s.t.}&u_j\le v_i+d_{i-1,j-1}, \quad\forall 0\leq i\leq m; 0\leq j\leq n
\end{array}
$$
and the following dual solution 
\begin{equation*}
\begin{array}{ll}
u_i=v_i=1-d_{i-1,n-1} & \mbox{for $i=0,...,m$} \\
u_j = 1-d_{m-1,n-1}+d_{m-1,j-1} &  \mbox{for $j=m+1,...,n$}.
\end{array}
\label{dualsolution}
\end{equation*}
This solution is dual-feasible: for $i,j\le m$ the triangle inequality gives
$$u_j-v_i=d_{i-1,n-1} -d_{j-1,n-1}\leq d_{i-1,j-1},$$
 whereas when $i\leq m<j$ the convex quadrangle inequality yields
$$u_j-v_i =  d_{m-1,j-1}-d_{m-1,n-1}+d_{i-1,n-1} \le d_{i-1,j-1}.$$
Replacing this solution in the dual objective function we obtain $D_{m,n}(\pi^n)$ once again, so that
strong duality implies that $(u,v)$ is dual optimal, and $z$ is a primal optimal transport.
Therefore $d_{m,n}=D_{m,n}(\pi^n)$.
\end{proof}

Exploiting the previous result, and since  the convex quadrangle inequality is a consequence of the monotonicity condition 
\eqref{eq_monotone}, we are led to consider the following 
constrained sequential optimization approach in which we restrict the $\pi^n$'s to satisfy this monotonicity as well as $\pi_n^n\geq\frac{1}{2}$.
Note that fixing $\pi^m$ and the distances $d_{i,j}$ for $m,i,j\leq n-1$, the expression 
$D_{m,n}(\pi^m)$ is a linear function of $\pi^n$, so that the following is again a linearly constrained quadratic programming problem 
\begin{equation*}
\begin{array}{rl}
\mbox{\sc (ms$_n$)} \hspace{0.5cm} \displaystyle{\min_{\pi^n\in \Delta^{\!n}}} &\pi^n_0+\displaystyle{\sum_{k=1}^n}\pi^n_k\, D_{k-1,n}(\pi^n) \\[2ex]
\mbox{s.t.} & \pi^n_k \le \pi^{n-1}_k,\qquad \forall k=0,...,n-1 \\
                  & \pi^n_n\ge \frac{1}{2} 
\end{array}
\end{equation*}

\subsection{Comparison of optimization strategies}\label{variedades}
We observe that {\sc (ms$_n$)} is a restricted version of {\sc (s$_n$)}, which is in turn a 
restricted version of {\sc (fh$_n$)}, so that their corresponding optimal values satisfy
$$\mbox{val}\mbox{\sc (fh$_n$)}\le \mbox{val}\mbox{\sc (s$_n$)} \le \mbox{val}\mbox{\sc (ms$_n$)}.$$
For $n=1$ all three  problems  share the same 
solution $\pi^{1}=(\frac{1}{2},\frac{1}{2})$ with optimal value $\frac{3}{4}$.
For $n=2$ the exact solution of {\sc (fh$_2$)} is $\pi^1 = (\sqrt{6}-2,3-\sqrt{6})$, 
$\pi^2 = (3\sqrt{6}-7,5-2\sqrt{6},3-\sqrt{6})$ with optimal value $30-12\sqrt{6}\sim 0.60612$. 
Note that the first stage optimal parameters $\pi^1$ in {\sc (fh$_2$)}  are different from those computed 
in {\sc (fh$_1$)}. This shows that the fixed horizon and sequential approaches are different.
As a matter of fact, the optimal solution for {\sc (s$_2$)} and  {\sc (ms$_2$)} is $\pi^2=(\frac{5}{14},\frac{1}{14},\frac{8}{14})$ 
with optimal value $17/28\sim 0.60714$. 

For $n\geq 3$  the analytic solutions  become increasingly harder to compute, so we proceed to compare
them numerically. Figure \ref{NLQPQPsimple} shows the reciprocal  $1/R_n$ of the optimal
residual bounds achieved by the three methods, plotted against $n$. The plot suggests that all
three methods yield algorithms with a rate of convergence $R_n\approx O(\frac{1}{n})$,
with slightly different slopes. Note that a higher slope 
corresponds to a faster convergence rate.  
Interestingly,  the numerical solutions reported by {\sc (s$_n$)} and {\sc (ms$_n$)} coincide
up to $n=100$, motivating the following
\begin{conjecture}\label{ConjectureMonotone}
The optimal sequence $\pi^n$ computed by {\sc (s$_n$)} satisfies  $\pi^n_n\geq\frac{1}{2}$ and $\pi_i^n\leq\pi_i^{n-1}$ for $i=0,\ldots,n-1$, so it coincides with the 
optimal sequence given by {\sc (ms$_n$)} with $\mbox{\rm val}\mbox{\sc (s$_n$)}= \mbox{\rm val}\mbox{\sc (ms$_n$)}$.
\end{conjecture}
 
\begin{figure}[t]
    \centering
    \includegraphics[scale=0.6]{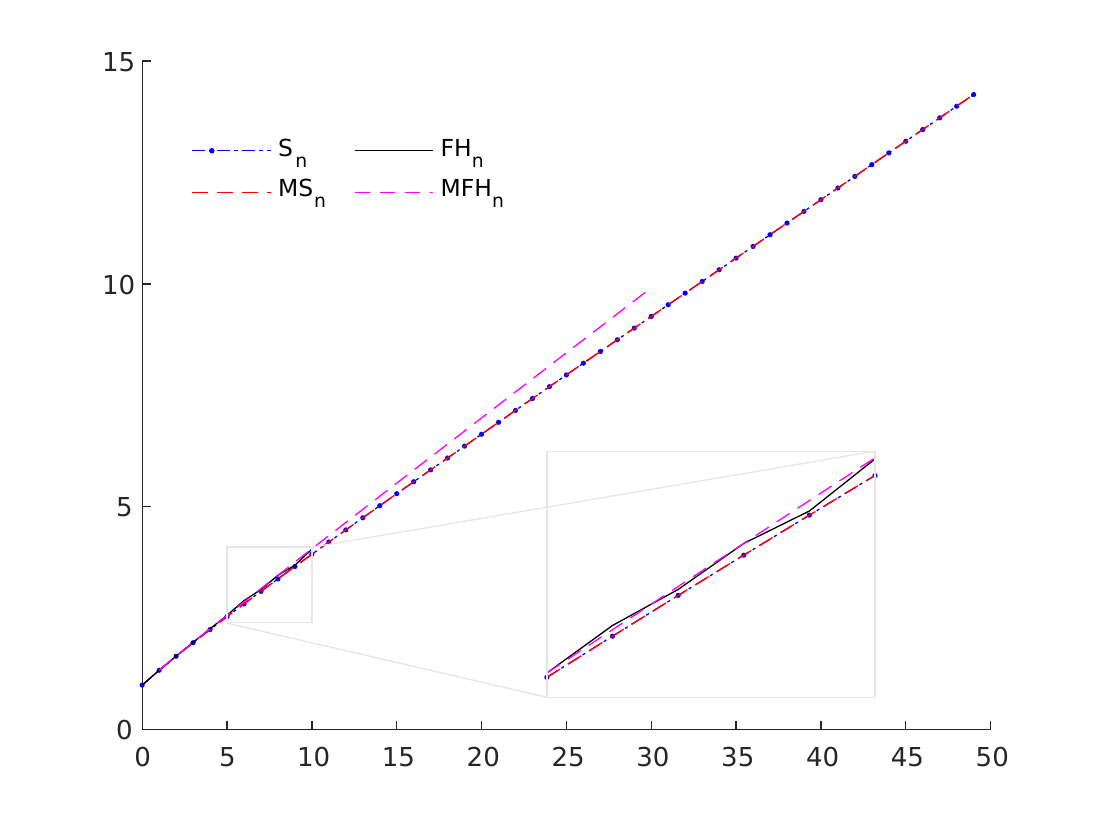}
    \caption{Plot of $1/R_n$ versus $n$ for {\sc (fh$_n$)}, {\sc (s$_n$)}, and {\sc (ms$_n$)}. 
    The zoomed image shows a numerical instability for {\sc (fh$_n$)} at $n=7$, possibly being trapped at a local minimum.}
    \label{NLQPQPsimple}
\end{figure}

Figure \ref{NLQPQPsimple} also plots the reciprocals of the optimal residual bounds achieved by a monotonic variant {\sc (mfh$_n$)} of the fixed horizon problem {\sc (fh$_n$)}
adding \eqref{eq_monotone} and \eqref{eqm2} as extra constraints, which together imply the structure  \eqref{zstructure} for the optimal transports $z^{k,m}$.
Empirically we observed that the solutions reported by {\sc (fh$_n$)} and {\sc (mfh$_n$)} coincide up to $n=6$. Beyond this limit {\sc (fh$_n$)} becomes hard to solve
and we could not check further if these problems yield the same solution as in the case of  {\sc (s$_n$)} and {\sc (ms$_n$)}.  
Proposing a conjecture based on such limited evidence would be hasty.

The small differences in the slopes in Figure \ref{NLQPQPsimple} suggest that 
the simpler sequential approaches  {\sc (s$_n$)}  and {\sc (ms$_n$)}  do not entail 
a significant loss with respect to the the fixed horizon scheme {\sc (fh$_n$)}. 
Moreover, this slight reduction in the proportionality constant is counterbalanced  
by the fact that the optimal coefficients $\pi^n$ can be computed more efficiently and do not need 
to be recalculated when the horizon $n$ changes. This provides a practical justification 
for using the sequential approach. 

Moreover, the problems {\sc (s$_n$)} and {\sc (ms$_n$)} are solved very quickly and to high accuracy, suggesting the existence of 
some further structure that could be exploited to derive the analytic properties of the solutions. However, a first 
inspection of these solutions did not reveal any clear and simple pattern. In fact, up to 
$n=10$ the optimal $\pi^{n}$ assigns a positive mass to all its components (see Figure \ref{plotpi}).
As expected, the mass concentrates mostly on the previous iterate $Tx^{n-1}$ but, most remarkably,  
the optimal solution also assigns significant mass to the initial point $y^0$, which decreases
as $n$ grows. 

\begin{figure}[t]
    \hspace*{-1.5cm}\includegraphics[scale=0.36]{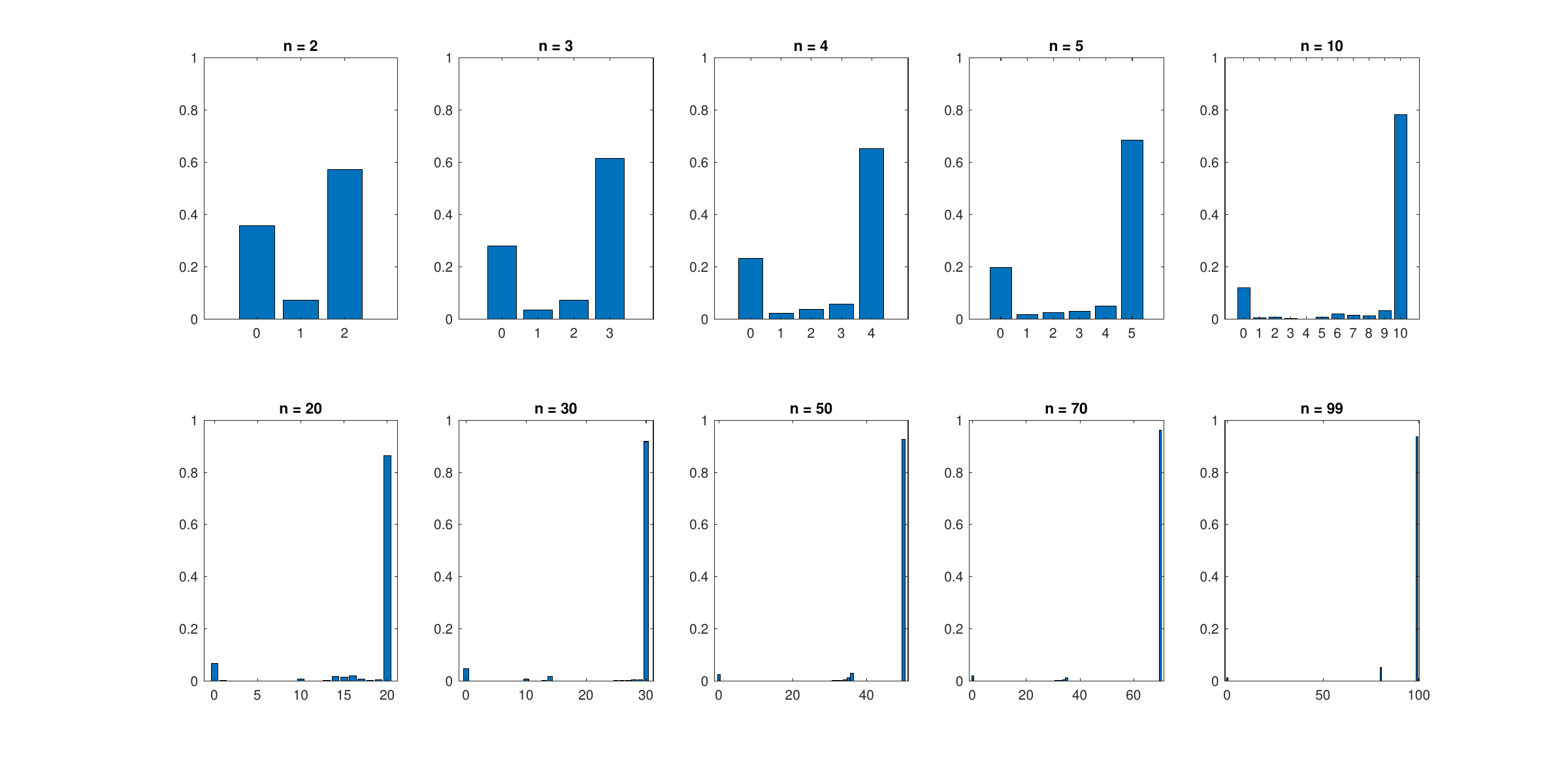}
    \caption{Optimal mass distribution with {\sc (ms$_n$)} for different values of $n$.}
    \label{plotpi}
\end{figure}

\subsection{Numerical optimization of particular Mann iterations} \label{specialMann}
In order to further understand the solutions given by the monotone sequential scheme {\sc (ms$_n$)}, 
we proceed to optimize the parameters $\pi^n$ for various particular instances of the general
Mann iterations. 
We already mentioned that by imposing additional structure on the coefficients $\pi^n$
one can recover the classical iterations of Krasnosel’ski{\u\i}, Halpern and Ishikawa. 
Also, motivated by Nesterov's acceleration schemes in convex optimization, we consider 
some variants that include additional terms.

Let $\{\alpha_k\}_{k\geq 1}, \{\beta_k\}_{k\geq 1}$ be two sequences in $[0,1]$, and denote by $\delta^n$ the Dirac mass at $n\in\NN$, i.e. $\delta^0 = (1,0,0,\ldots)$, $\delta^1 = (0,1,0,\ldots)$, and so on. We consider the following iterations which are particular cases of general Mann iterations.

\vspace{1ex}

\begin{itemize}
\item[{\sc (h)}]  {\em Halpern}:
$\pi^n =(1\!-\!\beta_n)\delta^0+\beta_n\delta^n$
$$x^{n} = (1-\beta_n)y^0+\beta_nTx^{n-1}$$

\item[{\sc (km)}] {\em Krasnosel’ski{\u\i}-Mann}:
$\pi^n = (1\!-\!\alpha_n)\pi^{n-1} + \alpha_n\delta^n$
$$x^n = (1-\alpha_n)x^{n-1}+\alpha_n Tx^{n-1}$$

\item[{\sc (th)}] {\em Twofold Halpern}:
$\pi^n=(1\!-\!\alpha_n\!-\!\beta_n)\delta^0+\beta_n\delta^{n-1}+\alpha_n\delta^n$
$$x^{n} = (1\!-\!\alpha_n\!-\!\beta_n)y^0+\beta_nTx^{n-2}+\alpha_n Tx^{n-1}$$

\item[{\sc (tkm)}] {\em Twofold Krasnosel’ski{\u\i}-Mann}:
$\pi^n = (1\!-\!\alpha_n\!-\!\beta_n)\pi^{n-1}+\beta_n\delta^{n-1}+\alpha_n\delta^n$
$$x^n = (1\!-\!\alpha_n\!-\!\beta_n)x^{n-1}+\beta_nTx^{n-2}+\alpha_nTx^{n-1}$$

\item[{\sc (kmh)}] {\em Krasnosel’ski{\u\i}-Mann-Halpern}:
$\pi^n = (1\!-\!\alpha_n\!-\!\beta_n)\delta^0+\beta_n\pi^{n-1}+\alpha_n\delta^n$
$$x^n = (1\!-\!\alpha_n\!-\!\beta_n) y^0+\beta_n x^{n-1}+\alpha_nTx^{n-1}$$

\item[{\sc (ekm)}] {\em Extra Krasnosel’ski{\u\i}-Mann}:
$\pi^n = (1\!-\!\alpha_n\!-\!\beta_n)\pi^{n-2} +\beta_n\pi^{n-1}+\alpha_n\delta^n$
$$x^n =(1\!-\!\alpha_n\!-\!\beta_n)x^{n-2}+\beta_nx^{n-1}+\alpha_nTx^{n-1}.$$
\end{itemize}

\vspace{2ex}
\noindent
Notice that the classical 2-step Ishikawa iteration 
$$\left\{\begin{array}{rcl}
x^{2n+1}&=&\!(1\!-\!\beta_n)x^{2n}\!+\beta_n Tx^{2n}\\
x^{2n+2}&=&\!(1\!-\!\alpha_n)x^{2n}\!+\alpha_n Tx^{2n+1}
\end{array}\right.,\quad 0\leq\alpha_n\leq\beta_n\leq 1.\leqno \mbox{\sc (ish)}
$$
is also a special case of  {\sc (ekm)}.  Moreover, for all the iterations above the coefficients $\pi^n$ are  
convex combinations of vectors in the unit simplex and only the entries from 0 to $n$ can take positive values.

Figure \ref{plotvariedades} plots the reciprocals $1/R_n$ of the residual bounds obtained
for all these iterations, with coefficients $\alpha_n$ and $\beta_n$ optimized in each
case by solving {\sc (ms$_n$)} with the corresponding additional structure of the $\pi^n$'s.
 For reference we also plot the curve 
{\sc (ms$_n$)} from Figure \ref{NLQPQPsimple}, which attains the best performance 
since there are no extra constraints on the coefficients. 
\begin{figure}[t]
    \centering
	\includegraphics[scale=0.8]{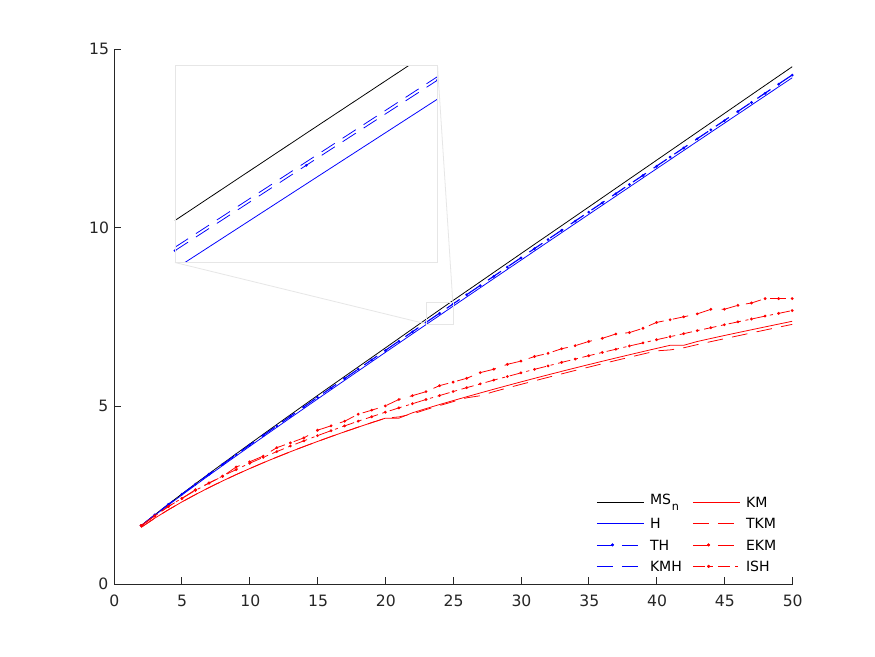}
    \caption{Plot of $1/\mbox{val}\mbox{\sc (ms$_n$)}$ vs. $n$ for particular instances of Mann's iterations. }
    \label{plotvariedades}
\end{figure}
The optimization was performed using Knitro with default options and then refined using Baron 
with the maximum number of iterations limited to 1000. 

We observe that the optimal {\sc (km)} does not attain the rate $O(1/n)$. 
As a matter of fact, it is known that  {\sc (km)} with constant stepsizes $\alpha_n\equiv\alpha$ 
attains the rate $O(1/\sqrt{n})$ ({\em cf.} Baillon \& Bruck \cite{baillon1996rate}). In the next section we   
prove that, even after optimizing the coefficients $\alpha_n$, {\sc (km)} cannot converge faster than $O(1/\sqrt{n})$.
Notably, these numerical results suggest that incorporating extra 
terms ---whether previous iterates  as in {\sc (ekm)} or images of previous iterates as in {\sc (tkm)}--- does not seem to
improve this rate.

 On the other hand, all the Halpern-type
methods ({\em cf.} the blue curves), which explicitly incorporate the starting point in the averaging, attain the faster rate $O(1/n)$. 
Furthermore, the slopes of the blue curves do not deviate much from the best value given by 
{\sc (ms$_n$)}. For instance, the average slope of the simplest
Halpern iteration is 0.2573, while for {\sc (ms$_n$)} is 0.2626, with a mere 2\% of relative difference.

Naturally, the optimal residuals decrease as we include more terms in the averaging process. 
However, the numerical results suggest that the performance of the simplest Halpern iteration that only takes $y^0$ and $Tx^{n-1}$
is comparable to the best general Mann iteration which involves all previous iterates:  in both cases we observe a rate $O(1/n)$ with only  some 
marginal gain in the proportionality constant which drops from $R_n\sim \frac{4}{n}$ to $R_n\sim \frac{10}{3n}$. In the next section we will formally prove that the rate $O(1/n)$ is the best one can expect from a 
general Mann iteration, whereas in Section \ref{optHalpern} we will discuss several variants of Halpern that attain this optimal 
rate, determining analytically the optimal coefficients $\beta_n$.

\section{Lower bounds on the convergence rates}\label{LB} 
The numerical computations suggest several conjectures. 
In particular, the Krasnosel'ski\u{\i}-Mann iterations and their variants
including extra terms seem to be systematically slower than the Halpern-like iterations.
Specifically, the optimal {\sc (km)} seems to attain a convergence rate of order $O(1/\sqrt{n})$,
whereas the optimal {\sc (h)} behaves as $O(1/n)$.  

In this section we formally prove that the 
fastest possible rate for general Mann iterations is in fact limited to $O(1/n)$.
It then follows that Halpern's iteration attains the optimal rate. 
We also prove that the best rate that can be achieved with an optimal 
Krasnosel'ski\u{\i}-Mann iteration is not better than $O(1/\sqrt{n})$, which is in fact
reached with constant stepsizes $\alpha_n\equiv\alpha$.

We stress that these lower bounds concern the worst-case for general normed spaces,
and faster  rates can be obtained by imposing further structure on the space $X$
({\em e.g.} Hilbert) or the map $T$ ({\em e.g.} a strict contraction).

\subsection{A lower bound for general Mann iterations}\label{lowerbound}
We begin by presenting a simple example showing that the best 
possible rate  for general Mann iterations in normed spaces is bounded from 
below by $O(1/n)$, independent of the averaging sequence $\pi^n$. 
The example map is a simple linear operator so that, even under this additional 
structure, one cannot expect a faster rate. Interestingly, the same map 
was recently used by Colao \& Marino \cite{colao2021} to establish lower bounds on $\|x^n-x^*\|$ for 
Halpern's iteration in $q$-uniformly smooth Banach spaces.

Let $T:\ell^\infty(\NN)\mapsto \ell^\infty(\NN)$ be the right-shift linear operator 
 given by 
\begin{equation}
T(x_0,x_1,x_2,\dots) = (0,x_0,x_1,x_2,\dots).
\label{Tshift1}
\end{equation} 
This map has a unique fixed point at $x^*=(0,0,0,\ldots)$, and the convex 
set $C=[0,1]^\NN$ is invariant under $T$ with $\diam(C)=1$. 
\begin{proposition}\label{LowerBoundMann} Consider a general Mann iteration \eqref{generalmann} for an arbitrary
sequence $(\pi^n)_{n\in\NN}$. Then,  the iterates
for the right-shift linear operator \eqref{Tshift1} started from $x^0=y^0=(1,1,1,\dots)$
satisfy 
$$\|x^n-Tx^n\|_\infty\geq\mbox{$\frac{1}{n+1}$}.$$
\end{proposition}
\begin{proof} Inductively one can check that the $n$-th  
iterate $x^n$ in every Mann iteration satisfies $x^n_i = 1$ for $i\geq n$, and therefore
\begin{equation*}
\|x^n-Tx^n\|_\infty \geq \varphi((x^n_i)_{i=0}^{n-1}) 
\label{Tshift1res}
\end{equation*}
where $\varphi:\RR^{n}\mapsto \mathbb{R}$ is the convex function 
$$\varphi(x_0,\ldots,x_{n-1}) = \max\{x_0,x_1-x_0,x_2-x_1,\ldots,x_{n-1}-x_{n-2},1-x_{n-1}\}.$$
The result will follow by showing that $\varphi$ attains its minimum at $\bar x_i =\frac{i+1}{n+1}$ for $i=0,\ldots,n-1$
with optimal value $\varphi(\bar x)=\frac{1}{n+1}$. Indeed, the coordinates of $\bar x$ are equidistant 
and all the terms in the maximum that define $\varphi(\bar x)$ are equal to $\frac{1}{n+1}$. 
Thus, denoting $e^k\in\RR^n$ for $k=0,\ldots,n-1$ the vector with $e^k_i=1$
if $i=k$ and $e^k_i=0$ otherwise, the subdifferential of $\varphi$ at $\bar x$ is 
$$\partial\varphi(\bar x)=\mbox{\rm co}\{e^0,e^1-e^0,e^2-e^1,\ldots,e^{n-1}-e^{n-2},-e^{n-1}\}.$$
Clearly $0\in \partial \varphi(\bar x)$ so that $\bar x$ is a minimizer with $\varphi(\bar x) = \frac{1}{n+1}$.   
\end{proof}

Note that while in this example the residual $\|x^n-Tx^n\|_\infty$ might converge to 0, the distance to the fixed point
$x^*$ is constant $\|x^n-x^*\|_\infty\equiv 1$. This illustrates the well-known fact that Mann iterates might not converge in norm to a fixed point, and the best one can expect in general is weak convergence, even when restricting to linear maps. 
However, 
strong convergence can be  guaranteed when $X$ is either a Hilbert space  ({\em cf.} \cite[Wittmann]{wittmann1992approximation})
 or a uniformly smooth Banach space ({\em cf.} \cite[Xu]{xu2002iterative}).

Let us also observe that although the lower bound and the numerical experiments agree on the order of
convergence $O(1/n)$, there is a gap between the lower bound $\|x^n-Tx^n\|_\infty\geq\frac{1}{n+1}$ and the upper
bounds which behave as $\frac{c}{n+1}$ with $c\approx 4$. In Section \ref{optHalpern} we will show that the optimal Halpern iteration 
achieves the latter bound with $c=4$, and that this constant is tight.

\subsection{A lower bound for Krasnosel'ski\u{\i}-Mann iterations}\label{KMlowerbound}
Our next result shows that for the general Krasnosel'ski\u{\i}-Mann iteration the best rate that can be achieved
is $O(1/\sqrt{n})$. As mentioned in the introduction, this optimal rate is in fact attained with constant stepsizes 
$\alpha_n\equiv\alpha$ (see Baillon \& Bruck \cite{baillon1996rate}). The following example considers again the right-shift 
operator, this time as a map acting on the space
$(\ell^1(\NN),\|\cdot\|_1)$.

\begin{proposition}
\label{lbkm}
Consider {\sc (km)} with arbitrary stepsizes $\{\alpha_n\}_{n\ge 0}$. Then, the iterates for the
right-shift operator \eqref{Tshift1} started from $x^0=(1,0,0,...)$  satisfy $\|x^n-Tx^n\|_1 \ge\frac{1}{\sqrt{n+1}}.$
\end{proposition}
\begin{proof}
We relegate the proof to the Appendix \ref{KMlowerboundProof} as it exploits the properties of binomial and Poisson binomial random variables, 
which are not relevant for the rest of this paper.
\end{proof}

\vspace{1ex}\noindent {\sc Remark 1.}
While the lower bound $\|x^n-Tx^n\| \ge\frac{1}{\sqrt{n+1}}$   
is the worst case in general Banach spaces, for Hilbert spaces it was observed by Baillon \& Bruck \cite[Section 9.4]{baillon1996rate} that the 
Krasnosel'ski\u{\i}-Mann iteration with constant stepsizes attains the faster
rate $\|x^n-Tx^n\|\sim o(\frac{1}{\sqrt{n}})$. Determining the exact rate in Hilbert spaces is still an open problem.

\section{Optimal Halpern iteration in normed spaces}\label{optHalpern}

Motivated by the lower bound in Section \ref{lowerbound}, and the numerical results in Section \ref{variedades} and Section \ref{specialMann},
we next focus on 
Halpern's iteration for which we determine analytically the optimal sequence of averaging parameters.
In particular, Theorem \ref{ThmHalpern}  proves that the asymptotic slope for this optimal Halpern iteration 
in Figure \ref{plotvariedades} converges exactly to 1/4 as $n\to\infty$.

Recall that Halpern's iteration is given by
\begin{equation}
x^{n} = (1-\beta_n)y^0 + \beta_nTx^{n-1}
\label{Halperndef}
\end{equation}
which corresponds to the averaging sequence $\pi^n=(1-\beta_n)\delta^0+\beta_n\delta^n$.
As shown next, the simple structure of the $\pi^n$'s implies that the minimizers of $\Psi_n(\pi)$ can be computed by solving either the
 fixed horizon problems {\sc (fh$_n$)} or the sequential schemes  
{\sc (s$_n$)} and {\sc (ms$_n$)}, which share the 
same optimal solution. In fact, the optimal $\beta_n$'s can be computed 
explicitly through a simple recursion.

\begin{theorem}\label{ThmHalpern}
Consider Halpern's iteration \eqref{Halperndef} for an arbitrary sequence $\beta_n\in [0,1]$. 
Then, the three optimization strategies {\sc (fh$_n$)}, {\sc (s$_n$)}, and {\sc (ms$_n$)} have the same optimal solution. More explicitly,
\begin{itemize}
\item[a)] the optimal $\beta_n\!$'s are given recursively by $\hat{\beta}_{n+1}\!\!=\!\frac{1}{2}(1\!+\!\hat{\beta}_n^2)$ with $\hat{\beta}_0\!=\!0$,\vspace{0.5ex}
\item[b)] the sequence $\hat{\beta}_n$ is increasing and the corresponding optimal bounds 
satisfy the recursion $\Ropt_{n+1}=\Ropt_{n}-\frac{1}{4}\Ropt_{n}^2$,\vspace{0.5ex}
\item[c)] the bounds $\Ropt_n$ are tight: there is a non-expansive $\hat{T}$ and a corresponding Halpern sequence such that $\|x^n-\hat{T}x^n\|=\Ropt_n$ for all
$n\in\NN$,\vspace{0.5ex}
\item[d)] the bounds satisfy $\Ropt_n\leq\frac{4}{n+4}$ with $\lim_{n\to\infty}(n+4)\Ropt_n=4$.
\end{itemize}

\end{theorem}

\begin{proof} In view of the simple structure $\pi_0^n=1-\beta_n$ and $\pi_n^n=\beta_n$ we have
 \begin{equation}\label{aux01}
 \mbox{$R_n=\sum_{i=0}^n\pi_i^nd_{i-1,n}=(1-\beta_n)+\beta_nd_{n-1,n}.$}
 \end{equation}
 On the other hand, for each $m\leq n$ there is a unique simple optimal transport.
 Namely, when $\beta_m\leq\beta_n$ the optimal transport is 
 $$\left\{\begin{array}{ccl}
 z_{0,0}&=&1-\beta_n\\
 z_{0,n}&=&\beta_n-\beta_m\\
 z_{m,n}&=&\beta_m
 \end{array}\right.
 $$
 all the other flows being null. Symmetrically, when $\beta_m\geq\beta_n$ the solution is
 $$\left\{\begin{array}{ccl}
 z_{0,0}&=&1-\beta_m\\
 z_{m,0}&=&\beta_m-\beta_n\\
 z_{m,n}&=&\beta_n
 \end{array}\right.
 $$
 so that both cases combined yield the recursive formula
$$ d_{m,n}=\vert \beta_m-\beta_n\vert +\min\{\beta_m,\beta_n\}d_{m-1,n-1}.$$
In particular
$$ d_{n-1,n}=\vert \beta_{n-1}-\beta_n\vert +\min\{\beta_{n-1},\beta_n\}d_{n-2,n-1}$$
  which plugged into \eqref{aux01} yields
$$R_n=\left\{
 \begin{array}{ll}
 \mbox{$1-\beta_n+\beta_n\left(\beta_n-\beta_{n-1}+\beta_{n-1} d_{n-2,n-1}\right)$}&\mbox{if }\beta_n\geq \beta_{n-1},\\
 \mbox{$1-\beta_n+\beta_n\left(\beta_{n-1}-\beta_n+\beta_n d_{n-2,n-1}\right)$}&\mbox{if }\beta_n\leq \beta_{n-1}.
 \end{array}\right.
$$ 
Now, since $d_{n-2,n-1}\leq 1$ and $\beta_{n-1}\leq 1$, the previous expression in the region $\beta_n\leq\beta_{n-1}$
is decreasing with respect to $\beta_n$,
and therefore the minimum of $R_n$ is achieved with $\beta_n\geq\beta_{n-1}$. 

Restricting to this latter case  and using \eqref{aux01} again, it follows that 
 \begin{eqnarray}\nonumber\label{recHal}
 R_n&=&1-\beta_n+\beta_n(\beta_n-\beta_{n-1}+\beta_{n-1}d_{n-2,n-1})\\
&=&1-\beta_n+\beta_n(\beta_n+R_{n-1}-1)\\\nonumber
&=&(1-\beta_n)^2+\beta_nR_{n-1}
 \end{eqnarray}
and then a simple induction yields the explicit formula
$$R_n=\mbox{$1-\sum_{k=1}^n(1-\beta_k)\prod_{i=k}^n\beta_i$}.$$
This expression could be plugged into the fixed horizon problems
{\sc (fh$_n$)}  and the sequential problem {\sc (s$_n$)} in order to 
find the optimal $\beta_n$'s. However, in this case
the problems admit a simple analytic solution. Indeed, observing 
the recursive structure \eqref{recHal} and noting that $R_{n-1}$ depends only
on the previous parameters $\beta_1,\ldots,\beta_{n-1}$, it follows that in the 
case of Halpern the fixed horizon and sequential approaches provide exactly
the same solution. Moreover, the optimal parameters can be found recursively 
by solving the trivial 1-dimensional quadratic problems
$$R_n=\min_{\beta_n\in [0,1]}(1-\beta_n)^2+\beta_nR_{n-1}$$
whose solution $\hat{\beta}_n$ satisfies $2(1-\beta_n)=R_{n-1}$. This yields the recursion 
\begin{eqnarray*}
2(1-\hat{\beta}_{n+1})&=&R_n\\
&=&(1-\hat{\beta}_n)^2+\hat{\beta}_nR_{n-1}\\
&=&(1-\hat{\beta}_n)^2+\hat{\beta}_n\,2(1-\hat{\beta}_n)\\
&=&1-\hat{\beta}_n^2
\end{eqnarray*}
which we rewrite in the form
\begin{eqnarray*}
\hat{\beta}_{n+1}&=&\mbox{$\frac{1}{2}(1+\hat{\beta}_n^2)$},\\
\Ropt_{n+1}&=&\mbox{$\Ropt_{n}-\frac{1}{4}\Ropt_{n}^2$}.
\end{eqnarray*}
This recursion clearly implies $\hat{\beta}_{n+1}\geq \hat{\beta}_{n}$ so that the fixed horizon and 
sequential approaches automatically satisfy the monotonicity condition \eqref{eq_monotone}.

The previous arguments show that for Halpern's iteration the optimal solutions 
of {\sc (fh$_n$)}, {\sc (s$_n$)}, and {\sc (ms$_n$)} coincide and minimize 
$\Psi_n(\pi)=R_n(\pi)$. These arguments also prove a) and b), while c) is a direct
consequence of Theorem \ref{TightMann}, so that it remains to prove d).

Defining $z_n = \frac{1}{4}\Ropt_n$ and using b) we get the recurrence $z_{n} = z_{n-1}(1-z_{n-1})$ with $z_0 = \frac{1}{4}$.
It follows that $z_n$ is decreasing and therefore $z_n<1$ for all $n\ge 0$. Rewriting this recurrence as
\begin{equation}\label{ttt}
\frac{1}{z_n}=\frac{1}{z_{n-1}}+\frac{1}{1-z_{n-1}}
\end{equation}
and using the inequality $\frac{1}{1-z_{n-1}}\geq 1$ we get  $\frac{1}{z_n}\geq\frac{1}{z_{0}}+n=4+n$
which yields precisely $\Ropt_{n}\le \frac{4}{n+4}$. In particular $z_n\searrow 0$
and the last term $\frac{1}{1-z_{n-1}}$ in \eqref{ttt} converges to 1, from which it follows the asymptotic
$\lim_{n\to\infty} n\, z_n=1$, and in turn $\lim_{n\to\infty}(n+4)\Ropt_n=4$. This completes the proof.
\end{proof}

Theorem \ref{ThmHalpern} combined with the global lower bound for Mann's iterations in Proposition \ref{LowerBoundMann},  
show that  Halpern's iteration with the recursive optimal stepsizes  $\hat{\beta}_{n+1}=\frac{1}{2}(1+\hat{\beta}_n^2)$ 
attains the optimal rate $O(1/n)$. 
The next result gives alternative conditions on the $\beta_n$'s that also guarantees the rate $O(1/n)$.
\begin{proposition}\label{propHal2}
Consider Halpern's iteration with $\beta_n\in[0,1]$, and assume that for some constants $a,\kappa$ such that 
$1\le a+1\le \kappa$ and $\kappa\geq 4$, we have
\begin{equation}\label{eq_propHal2}
\mbox{$(1-\beta_n)^2+\frac{\kappa}{n+a}\beta_n \le \frac{\kappa}{n+a+1}$.}
\end{equation}
 Then $\|x^n-Tx^n\| \le \frac{\kappa}{n+a+1}$ for all $n\ge 0$.
\end{proposition}
\begin{proof} We prove inductively that $R_n\leq\frac{\kappa}{n+a+1}$. For $n=0$ this holds trivially since $R_0=\diam(C)=1\leq \frac{\kappa}{a+1}$. 
Assuming $R_{n-1}\le \frac{\kappa}{n+a}$, equation \eqref{recHal} gives
$$R_n = (1-\beta_n)^2+\beta_nR_{n-1} \le \mbox{$(1-\beta_n)^2+\frac{\kappa}{n+a} \beta_n \le \frac{\kappa}{n+a+1}$}$$
completing the induction. 
\end{proof}

\vspace{1ex} \noindent {\sc Remark 2.}
The optimal rate $\|x^n-Tx^n\|\approx O(1/n)$ for Halpern's iteration was already obtained by 
Sabach \& Shtern \cite[Lemma 5]{sabach2017first} with stepsizes $\beta_n=\frac{n}{n+2}$
and with the explicit estimate $\|x^n-Tx^n\|\leq \frac{4}{n+1}$. 
Although their result was presented in $\mathbb{R}^d$, the proof is valid in any normed space. 
This result also follows from Proposition \ref{propHal2}
since \eqref{eq_propHal2} holds with $a=0$ and $\kappa=4$. 
On the other hand, we recall that for $\beta_n=\frac{n}{n+2}$ we have the 
tight bound $R_n=\frac{4}{n+1}(1-\frac{H_{n+2}}{n+2})$ with  $H_n=\sum_{k=1}^n\frac{1}{k}$ 
the $n$-th harmonic number (see Bravo {\em et al.} \cite{bcc2020}). 
This tight bound $R_n$ is very close to the optimal bound $\Ropt_n$ in Theorem \ref{ThmHalpern}. 
Numerically we observe that $R_n/\Ropt_n$ increases until $n=4$ attaining a maximal value of 1.05223, after which it decreases
converging asymptotically to 1. In other words $\beta_n = \frac{n}{n+2}$  and the optimal 
scheme $\hat{\beta}_n$ are asymptotically equivalent for  $n$ large. Moreover, both  $R_n$ and $\Ropt_n$ are tight for the corresponding $\beta_n$'s.\vspace{1ex}

\vspace{1ex} \noindent {\sc Remark 3.}
For $\kappa<4$ no sequence $\beta_n$ satisfies \eqref{eq_propHal2}.
With $\kappa = 4$ and $a=3$ the condition is equivalent to
$\vert\beta_n-\frac{n+1}{n+3}\vert \le \epsilon_n$ with $\epsilon_n=\frac{2}{(n+3)\sqrt{n+4}}$
and implies the bound $R_n\le \frac{4}{n+4}$.  
This holds in particular for $\beta_n = \frac{n+1}{n+3}$ which coincide with the stepsizes in Sabach \& Shtern \cite{sabach2017first} except that 
they are shifted by one and provide a slightly smaller  bound for the residuals. 
One can also prove that for $n\ge 1$ the optimal coefficients $\hat{\beta}_n$ satisfy $0\le \hat{\beta}_n-\frac{n+1}{n+3} \le \epsilon_n$ so
that we recover the bound in Theorem \ref{ThmHalpern}\,d).  \vspace{1ex}

\vspace{1ex} \noindent {\sc Remark 4.}
By rewriting the recursion $\hat{\beta}_{n+1}=\frac{1}{2}(1+\hat{\beta}_n^2)$ in terms of
the complementary values $\hat{\alpha}_n=1-\hat{\beta}_n$ we obtain $\hat{\alpha}_{n+1}=\hat{\alpha}_n-\frac{1}{2}\hat{\alpha}_n^2$.
As in the last argument in the proof of Theorem \ref{ThmHalpern} it follows that $n\,\hat{\alpha}_n\to 2$ as $n\to \infty$, and therefore 
$\hat{\alpha}_n\to 0$ with $\sum_{n}\hat{\alpha}_n=\infty$ and $\sum_{n}\vert \hat{\alpha}_{n+1}-\hat{\alpha}_n\vert <\infty$. Hence, from
Wittmann \cite{wittmann1992approximation} we conclude that when $X$ is a Hilbert space and $\Fix(T)\neq\phi$ the Halpern 
iterates $x^n$ converge in norm towards $P_{\Fix(T)}(y^0)$ the projection of $y^0$ 
onto the set of fixed points of $T$. Moreover, we also have $\hat{\alpha}_{n+1}/\hat{\alpha}_n\to 1$ so that 
Theorem 3.1 in Xu \cite{xu2002iterative} implies that strong convergence to a fixed point 
also holds when $X$ is a uniformly smooth Banach. A similar result can be derived from  \cite[Theorem 1]{reich1994} and  \cite[Remark 1]{reich1994}  
whenever $X$ has a weakly sequentially continuous duality map.

\subsection{Previous results for Halpern iteration in Hilbert spaces}\label{secPrevHalpern}

By restricting to Hilbert spaces, Lieder \cite{lieder2021convergence} recently established that Halpern 
iterates 
with stepsize $\beta_n = \frac{n}{n+1}$ achieve the accelerated rate with 
$$\|x^n-Tx^n\| \le \mbox{$\frac{2\|x^0-x^*\|}{n+1}$}, \qquad \forall n\ge 1.$$
Lieder also showed that this bound is sharp: for any fixed $n\in \mathbb{N}$ 
there exists a non-expansive map $T:\mathbb{R}\mapsto \mathbb{R}$, which depends on $n$, that attains the equality.
Lieder's bound improves by a factor 4 the proportionality constant in the bound in
Theorem \ref{ThmHalpern} d). The paper presents two proofs: a direct algebraic proof using 
an ad-hoc  weighted sum that provides a Lyapunov function, and a second proof using techniques of  Performance Estimation 
Problems (PEP). Both proofs strongly exploit the paralellogram identity. 

On the other hand, as shown in Bravo {\em et al.} \cite{bcc2020}, for $\beta_n = \frac{n}{n+1}$ the best bound that
can be achieved in general Banach spaces is $R_n=\frac{H_{n+1}}{n+1}\sim O(\frac{\ln n}{n})$
showing that the PEP approach is intrinsically restricted to a Hilbert setting. On the positive side, this shows that the bounds $R_n$ can be sharpened if we 
restrict the space $X$ on which the non-expansive map $T$ is defined.

Another accelerated iteration for finding zeros of co-coercive operators in Hilbert spaces was recently 
proposed by Kim \cite{kim2021accelerated}. Recall that $M\!:\!{\mathcal H}\!\to\!{\mathcal H}$ is co-coercive with parameter $\mu$ iff
$$\langle Mx-My,x-y\rangle\geq\mu \|Mx-My\|^2\quad\forall~x,y\in{\mathcal H}.$$
It is well known that  this is equivalent to $T\triangleq I-2\mu M$ being non-expansive. Conversely, 
a map $T:{\mathcal H}\to{\mathcal H}$ is non-expansive iff $M=I-T$ is $\frac{1}{2}$-cocoercive
so that, in a Hilbert setting, finding zeros of cocoercive operators is equivalent to finding fixed point for non-expansive maps. 
All of this is well known, and we just recall it for the reader's convenience. By considering the general iteration 
$$x^{n+1} = x^n -\mu\sum_{k=0}^n h_{k+1,n+1}Mx^{k},$$
where $h$ is a matrix of stepsize values, and after inspecting the optimal stepsizes
obtained using techniques based on PEP,  
Kim \cite{kim2021accelerated} proposed the following inertial method started from $x^0=y^0=x^{-1}$
\begin{equation}
\label{kimmethod}
\begin{array}{rcl}
y^{k+1} &=& (I-\mu M)x^k,\\
x^{k+1} &=& \mbox{$y^{k+1}+\frac{k}{k+2}(y^{k+1}-y^k)-\frac{k}{k+2}(y^k-x^{k-1}),$}
\end{array}
\end{equation}
and showed that $\|Mx^n\|\le \frac{2\|x^0-x^*\|}{n+1}$. 
While Kim's method looks different from Halpern, we show below that it is in fact equivalent to the 
iteration studied in Lieder \cite{lieder2021convergence}. This connection between inertial techniques 
and Halpern does not seem to have been noticed earlier, and might shed additional light on the mechanisms involved in the 
acceleration of fixed point iterations.

\begin{proposition}\label{propKim}
Let $T:C\mapsto C$ be non-expansive. Then the iteration 
\eqref{kimmethod} applied to the $\frac{1}{2}$-cocoercive operator $M = I-T$
coincides with the classical Halpern iteration $x^{n} = \frac{1}{n+1}x^0+\frac{n}{n+1}Tx^{n-1}$.
\label{kimequalHalpern}
\end{proposition}

\begin{proof}
Taking  $M=I-T$ and $\mu=1/2$ in the first equation of \eqref{kimmethod} gives
$$y^{k+1} = (I-\mu M)x^k = \mbox{$\frac{1}{2}$}(x^k+Tx^k)$$
which substituted into the definition of $x^{k+1}$ yields
\begin{eqnarray}
x^{k+1} &=& \mbox{$\frac{1}{2}(x^k+Tx^k)+\frac{k}{k+2}\left(\frac{1}{2}(x^k+Tx^k)-\frac{1}{2}(x^{k-1}+Tx^{k-1})\right)$}\nonumber \\
&& \mbox{$\qquad -\frac{k}{k+2}\left(\frac{1}{2}(x^{k-1}+Tx^{k-1})-x^{k-1})\right)$} \nonumber\\
&=& \mbox{$\frac{k+1}{k+2}(x^k+Tx^k)-\frac{k}{k+2}Tx^{k-1}$}.\nonumber
\end{eqnarray}
Multiplying by $k+2$ and rearranging we obtain
$$(k+2)x^{k+1}-(k+1)Tx^k = (k+1)x^k-kTx^{k-1}$$
which shows that the sequence $z^k=(k+1)x^k-kTx^{k-1}$ is constant. 
Hence $z^n=z^0$ which yields $(n+1)x^n-nTx^{n-1}=x^0$
and the result follows.
\end{proof}

\vspace{1ex} \noindent {\sc Remark 5.}
During the review process of this paper we learned that an equivalent result was obtained by Ryu \& Yin  (see \cite[Chapter 12.2, Theorem 18]{ryu2021large}). 
We thank an anonymous referee for pointing out this relevant reference.

\subsection{Optimal Halpern iteration for affine maps}
In this section, we investigate Halpern's iteration 
in the context of a linear non-expansive map $T:X\mapsto X$, with $X$ a general normed space. We will show that in this case we recover the optimal bound $\frac{2}{n+1}$ known for general non-expansive maps in Hilbert spaces. We remark that the  linear case also covers the setting of an affine operator $T$ with $\Fix(T)\neq\phi$, as we can conveniently translate the origin.

Let us fix a sequence $\beta = (\beta_n)_{n}$ with $\beta_0=0$, and denote 
\begin{equation*}\label{pikn}
\mbox{$\Pi^n_k(\beta) = \prod_{l=k}^n\beta_l$}
\end{equation*}
where by convention the product is set to $1$ for $k>n$. 
Using the linearity of $T$, a straightforward induction shows that the
Halpern iterates with initial data $x^0=y^0$ can be expressed as
\begin{equation}\label{Halpernlineal}
x^n = \sum_{k=0}^n \mbox{$(1-\beta_k)\Pi_{k+1}^n(\beta)$}\,T^{n-k}x^0,
\end{equation}
and by linearity of $T$
\begin{equation}\label{HalpernlinealT}
Tx^n = \sum_{k=0}^n \mbox{$(1-\beta_k)\Pi_{k+1}^n(\beta)$}\,T^{n-k+1}x^0.
\end{equation}

\begin{theorem}\label{theoremHalpernlineal}
Let  $T:X\mapsto X$ be linear non-expansive and $x^*\in\Fix(T)$. 
Let also $\Theta_n:[0,1]^n\mapsto \mathbb{R}$ be defined by 
\begin{equation}\label{theta}
\mbox{$\Theta_n(\beta) = 1\!-\!\beta_n+\Pi_1^n(\beta) + \sum_{k=1}^{n} \left\vert(2-\beta_{k-1})\beta_{k}-1\right\vert\, \Pi_{k+1}^n(\beta).$}
\end{equation}
Then, the Halpern iterates with initial data $x^0=y^0$ satisfy
\begin{equation}\label{Halpernlinealineq}
\|x^n-Tx^n\| \le \|x^0\!-x^*\|\,\Theta_n(\beta).
\end{equation}
Moreover, this bound is tight: for the right-shift operator $T:\ell^1(\NN)\mapsto \ell^1(\NN)$ 
$$T(x_0,x_1,x_2,\ldots)=(0,x_0,x_1,x_2,\ldots),$$
the iterates started from $x^0=y^0=(1,0,0,\ldots)$ satisfy \eqref{Halpernlinealineq}
with equality.
\end{theorem}

\begin{proof} From  \eqref{Halpernlineal} we have
\begin{eqnarray*}
x^n\!-\!Tx^n&=&\sum_{k=0}^n (1-\beta_k)\Pi^n_{k+1}(\beta)\,(T^{n-k}x^0\!-T^{n-k+1}x^0) \\
&=&\sum_{k=0}^n (1-\beta_k)\Pi^n_{k+1}(\beta)\,(T^{n-k}(x^0-x^*)-T^{n-k+1}(x^0-x^*)) \\
&=&\sum_{k=1}^n [(2-\beta_{k-1})\beta_{k}-1]\Pi^n_{k+1}(\beta)\,T^{n-k+1}(x^0-x^*) \\
&&{}+(1-\beta_n)(x^0-x^*)-\Pi^n_{1}(\beta)T^{n+1}(x^0-x^*).
\end{eqnarray*}
Here, the first line comes from the equations \eqref{Halpernlineal} and \eqref{HalpernlinealT}. In the second line we add $T^{n-k+1}x^*$ and subtract $T^{n-k}x^*$, both of which coincide with $x^*$. In the third line, we re-arrange the terms in the sum. Then, \eqref{Halpernlinealineq} follows by taking the norm and using the triangle inequality and non-expansivity of $T$. 

To prove the tightness of \eqref{Halpernlinealineq} let $\delta^k$ be the $k$-th canonical vector in $\ell^1(\NN)$. 
Since $x^0=\delta^0$ and $T$ is linear with $T\delta^k = \delta^{k+1}$, from \eqref{Halpernlineal} we get
$$x^n = \sum_{k=0}^{n} (1-\beta_k)\Pi^n_{k+1}(\beta)\,T^{n-k}\delta^0 = \sum_{k=0}^{n} (1-\beta_k)\Pi^n_{k+1}(\beta)\delta^{n-k},$$
and therefore
\begin{align*}
\|x^n-Tx^n\|_1 & =  \vert 1-\beta_n\vert  + \sum_{k=1}^n\vert (1\!-\!\beta_k)\Pi^n_{k+1}(\beta)-(1\!-\!\beta_{k-1})\Pi^n_{k}(\beta)\vert +\Pi^n_1(\beta)\\
 & = 1-\beta_n+\Pi^n_1(\beta) + \sum_{k=1}^n\vert (1-\beta_k)-(1-\beta_{k-1})\beta_k\vert \Pi^n_{k+1}(\beta)\\
&=\Theta_n(\beta)
\end{align*}
The result follows since $T$ has a fixed point at $x^*=0$ with $\|x^0-x^*\|_1 = 1$. 
\end{proof}

Since the bound \eqref{Halpernlinealineq} is always tight for any choice of the $\beta_n$'s, the best possible Halpern iteration 
for linear maps is obtained by minimizing the function $\Theta_n(\cdot)$. As shown below, this minimum can be computed 
explicitly.
\begin{proposition}\label{minTheta}
The minimum of the function $\beta\mapsto\Theta_n(\beta)$ defined by \eqref{theta} is attained
at $\beta^*_k = \frac{k}{k+1}$ for $k=1,...,n$ with value $\Theta(\beta^*)=\frac{2}{n+1}$. 
\end{proposition}
\begin{proof}
Let $y_k = \prod_{l=k}^n \beta_l$ for $k=1,...,n$, $y_{n+1} = 1$, and $y_0=0$. Note that we have the inverse relation $\beta_k = \frac{y_{k}}{y_{k+1}}$, for $k=0,...,n$.  
Consider the optimization problem
\begin{equation*}
(P_y) \;\; \min_{\begin{array}{c}y_{0},...,y_{n+1}\in [0,1]\\ y_0 = 0, y_{n+1}=1, \end{array}}1-y_n+y_1+\sum_{k=1}^{n}\vert y_{k+1}-2y_{k}+y_{k-1}\vert
\label{proby}
\end{equation*}
which is a relaxation of the original problem $\min_{\beta\in [0,1]^{n}}\Theta_n(\beta)$ as we are not considering the constraints $y_{k}\le y_{k+1}$ for $k=0,...,n$. 

Now let us make a second change of variables: $\omega_{k+1} = y_{k+1}-2y_{k}+y_{k-1}$ for $k=1,...,n$ and $w_1=y_1$.
This is a linear transformation that can be easily reverted as $y_k = \sum_{i=1}^k (k-i+1)\omega_i$ for $k=1,...,n+1$. Consequently,
$$1-y_n = y_{n+1}-y_n = \sum_{k=0}^{n+1}(n-k+2)\omega_k-\sum_{k=1}^n(n-k+1)\omega_k = \sum_{k=1}^{n+1}\omega_k.$$
Moreover, from the relation $y_{n+1}=1$ we can compute $\omega_1$ as
$$\mbox{$\omega_1 = \frac{1}{n+1}\left(1-\sum_{k=2}^{n+1}(n-k+2)\omega_k\right).$}$$
This change of variables transforms the objective function of $(P_y)$ into  
\begin{equation}
\label{probomega}
\mbox{$\frac{2}{n+1}+\sum_{k=2}^{n+1}\left[ \left(1-2\frac{n-k+2}{n+1}\right)\omega_k+\vert \omega_k\vert\right].$}
\end{equation}
The coefficients $1-2\frac{n-k+2}{n+1}$ belong to $(-1,1)$ for all $k=2,...,n+1$. Hence the terms in the sum 
\eqref{probomega} are always non-negative and this expression attains its minimum with $\omega_k = 0$ for $k=2,...,n+1$. 
Consequently $\omega_1 = \frac{1}{n+1}$ and the solution to $(P_y)$ is $y_{k} = \frac{k}{n+1}$ for $k=1,...,n+1$
which yield the coefficients $\beta^*_k=\frac{k}{k+1}$. 
This solution satisfies the constraint $y_k\le y_{k+1}$ so that the optimal value of the relaxation $(P_y)$ 
coincides with $\Theta(\beta^*)$ and $\beta^*$ is optimal for $\Theta(\cdot)$.
\end{proof}

As a direct consequence of Theorem \ref{theoremHalpernlineal} and Proposition \ref{minTheta},
we derive the following optimal Halpern iteration for linear maps.
\begin{theorem} The optimal Halpern iteration for linear maps is obtained by choosing the stepsizes 
$\beta_n =\frac{n}{n+1}$, which attains the tight bound
\begin{equation}\label{Halpernlinealboundclasical}
\|x^n-Tx^n\| \le\mbox{$ \frac{2\|x^0-x^*\|}{n+1}$},\qquad \forall n\ge 1.
\end{equation}
With this choice the iterates are given by
$$x^n=\mbox{$\frac{1}{n+1}x^0+\frac{n}{n+1}Tx^{n-1}=\frac{1}{n+1}\sum_{k=0}^nT^kx^0.$}$$
\end{theorem}

A remarkable and unexpected fact is that the optimal coefficients $\beta_n=\frac{n}{n+1}$ and the bound
\eqref{Halpernlinealboundclasical} coincide exactly with the coefficients and bound for 
general non-expansive maps in Hilbert spaces obtained in Lieder \cite{lieder2021convergence}. 
Here the space $X$ is not required to be Hilbert, but the 
map $T$ is assumed linear. Moreover, as illustrated by the next example, the bound \eqref{Halpernlinealboundclasical} remains tight in Hilbert spaces.

\vspace{1ex}
\noindent
{\sc Example.} Consider the map $T$ (which depends on $n$) given by the rotation 
in $(\mathbb{R}^2,\|\cdot\|_2)$ with angle $\theta_n = \frac{\pi}{n+1}$. This is clearly a 
linear map with a unique fixed point at the origin, and for all $k\ge 0$ we have
$$\mbox{$T^k x^0 = \left[\begin{matrix}
\cos k\theta_n & -\sin k\theta_n \\
\sin k\theta_n & \cos k\theta_n
\end{matrix}\!\right]x^0$.}$$
Halpern's iteration with $\beta_k =\frac{k}{k+1}$ started from $x^0 = \left(\begin{matrix} 1\\0\end{matrix}\right)$
satisfies
$$x^n = \frac{1}{n\!+\!1}\sum_{k=0}^n
\mbox{$\left(\begin{matrix}
\cos k\theta_n \\ \sin k\theta_n
\end{matrix} \!\right)$}$$ 
and a straightforward telescoping yields
$$
\|x^n-T x^n\|_2 = \left\|\frac{1}{n\!+\!1} \sum_{k=0}^n
\mbox{$\left(\left(\begin{matrix}
\cos k\theta_n \\ \sin k\theta_n
\end{matrix} \!\right)-\left(\begin{matrix}
\cos (k\!+\!1)\theta_n \\ \sin (k\!+\!1)\theta_n
\end{matrix} \!\right) \right)$}\right\|_2= \frac{2}{n+1}. $$

\begin{appendices}

\section{Tightness of optimal transport bounds}\label{SecTightMann}

In this section we present the proof of Theorem \ref{TightMann}, establishing the tightness 
of the optimal transport bounds in general Mann iterations, and therefore 
the equality $\Psi_n(\pi)=R_n(\pi)$.
\begin{proof} [Proof of Theorem \ref{TightMann}]
Let $z^{mn}$ and $u^{mn}$ be optimal solutions for $({\mathcal P}_{m,n})$ and $({\mathcal D}_{m,n})$.
Setting $u_{i}^{mn}={\displaystyle \min_{0\leq k\leq n}}u_{k}^{mn}+d_{k-1,i-1}$ for $i>n$,
and using the triangle inequality, we get
\begin{equation}\label{modulo}
\vert u_{i}^{mn}-u_{j}^{mn}\vert \leq d_{i-1,j-1}\mbox{ for all }i,j\in\NN.
\end{equation}
In particular all the $u_i^{mn}$'s are within a distance at most 1 and,
since the objective function in $({\mathcal D}_{m,n})$ is invariant by 
translation, we may further assume that $u_i^{mn}\in[0,1]$ for all $i\in\NN$.

Let $\D$ be the set of all pair of integers $(m,n)$ with $-1\leq m\leq n$, and consider the 
unit cube $C=[0,1]^\D$ in the space $(\ell^\infty(\D),\|\cdot\|_\infty)$.
For every integer $k\in\NN$ define $y^{k}\in C$ as
 \begin{equation}\label{yyy}  
 \forall (m,n)\in \D\qquad y^k_{m,n}=\left\{\begin{array}{cl}
 d_{k-1,n}&\mbox{if } -\!1=m\leq n\\
 u^{mn}_{k}&\mbox{if }\hspace{2.5ex} 0\leq m\leq n
 \end{array}\right.
\end{equation}
and  a corresponding sequence $x^{k}\in C$ given by
 \begin{equation}
 \label{rec} \mbox{$ x^k=\sum_{i=0}^k\pi_i^ky^{i}$}.
 \end{equation}
 
 We claim that $\|y^{m+1}-y^{n+1}\|_\infty\leq d_{m,n} =\|x^m-x^n\|_\infty$ for all $0\leq m\leq n$. 
 Indeed, using the triangle inequality and  \eqref{modulo}  we get
$$
\left\{\begin{array}{ll}
\vert y^{m+1}_{-1,n'}-y^{n+1}_{-1,n'}\vert =\vert d_{m,n'}-d_{n,n'}\vert \leq d_{m,n}&\mbox{~~~~if $-1=m'\leq n'$}\\[1ex]
\vert y^{m+1}_{m',n'}-y^{n+1}_{m',n'}\vert =\vert u^{m'n'}_{m+1}\!-u^{m'n'}_{n+1}\vert \leq d_{m,n}&\mbox{~~~~if\hspace{2ex} $0\leq m'\leq n'$}
\end{array}\right.
$$
which together imply
\begin{equation}\label{ydes}
\|y^{m+1}\!-y^{n+1}\|_\infty\leq d_{m,n}.
\end{equation}
Also, selecting an optimal transport $z^{mn}$ for $({\mathcal P}_{m,n})$ 
we have
\begin{eqnarray}\label{diferencia}\nonumber
x^m-x^n&=&\mbox{$\sum_{i=0}^{m}\pi_{i}^my^{i}-\sum_{j=0}^{n}\pi_j^ny^{j}$}\\
&=&\mbox{$\sum_{i=0}^{m}\sum_{j=0}^{n}z^{mn}_{i,j}(y^{i}-y^{j})$}
\end{eqnarray}
so that the triangle inequality and \eqref{ydes} yield
\begin{equation}\label{siete}
\mbox{$
\|x^m\!-x^n\|_\infty\leq 
\sum_{i=0}^{m}\sum_{j=0}^{n}z^{mn}_{i,j}d_{i-1,j-1}=d_{m,n}.
$}
\end{equation}
On the other hand, considering the $(m,n)$-coordinate in \eqref{diferencia}, the complementary slackness 
\eqref{cs} gives (recall that we are in the case $0\leq m\leq n$)
 \begin{eqnarray*}
 \vert x_{m,n}^m-x_{m,n}^n\vert &=&\vert \mbox{$\sum_{i=0}^{m}\sum_{j=0}^{n}z^{mn}_{i,j}(y_{m,n}^{i}-y_{m,n}^{j})$}\vert \\
 &=&\vert  \mbox{$\sum_{i=0}^{m}\sum_{j=0}^{n}z^{mn}_{i,j}(u^{mn}_i-u^{mn}_j)$}\vert \\
 &=&\mbox{$\sum_{i=0}^{m}\sum_{j=0}^{n}z^{mn}_{i,j}d_{i-1,j-1}=d_{m,n}$}
 \end{eqnarray*} 
 which combined with \eqref{siete} yields $\|x^m-x^n\|_\infty=d_{m,n}$ as claimed. 
 
Define $T:S\to C$ on the set $S=\{x^k:k\in\NN\}\subseteq C$ by $Tx^k=y^{k+1}$, so that $T$
is non-expansive. Since $\ell^\infty(\D)$ as well as the unit cube $C$ are hyperconvex, then by Theorem 4 in Aronszajn \& Panitchpakdi \cite{aronszajn1956extension}, $T$ can be extended to a non-expansive map $T:C\to C$ and then
\eqref{rec} is precisely a Mann sequence which attains all the bounds 
$\|x^m-x^n\|_\infty=d_{m,n}$ with equality.

It remains to prove that $\|x^n-Tx^n\|_\infty=R_n$.
The upper bound follows again using the triangle inequality and \eqref{ydes} since
$$\|x^n-Tx^n\|_\infty=\mbox{$\|\sum_{i=0}^n\pi^n_i(y^{i}-y^{n+1})\|_\infty$}\leq \mbox{$\sum_{i=0}^n\pi_i^nd_{i-1,n}=R_n$}.$$
For the reverse inequality, we look at the coordinate $(-1,n)$ so that
\begin{eqnarray*}
\|x^n-Tx^n\|_\infty&=&\mbox{$\|\sum_{i=0}^n\pi^n_iy^{i}-y^{n+1}\|_\infty$}\\
&\geq&\vert \mbox{$\sum_{i=0}^n\pi^n_iy_{-1,n}^{i}-y_{-1,n}^{n+1}\vert $}\\
&=&\vert \mbox{$\sum_{i=0}^n\pi^n_id_{i-1,n}-d_{n,n}\vert $}=R_n
\end{eqnarray*}
which completes the proof.
 \end{proof}

\vspace{1ex} \noindent {\sc Remark 6.}
As in Bravo {\em et al.}  \cite{bcc2020} we observe that the map $T$ must have a fixed point in $C$, and
also that it can be extended to the full space $\ell^\infty(\D)$.

\section{Lower bound for Krasnosel'ski\u{\i}-Mann iterations.}\label{KMlowerboundProof}

In this Appendix we prove Proposition \ref{lbkm} by showing a non-expansive linear operator $T$ for 
which the Krasnosel'ski\u{\i}-Mann sequence $x^{n+1} = (1-\alpha_n)x^n +\alpha_nTx^n$
satisfies $\|x^n-Tx^n\| \ge \frac{1}{\sqrt{n+1}}$, independently of the stepsizes 
$\{\alpha_n\}_{n\ge 0}$. 

\begin{proof} [Proof of Proposition \ref{lbkm}]
Let $T$ be the right-shift operator \eqref{Tshift1} considered as a map acting on $(\ell^1(\mathbb{N}),\|\cdot\|_1)$. 
Fix an arbitrary sequence of stepsizes $\alpha_n$ and consider the corresponding Krasnosel'ski\u{\i}-Mann sequence started from $x^0=(1,0,0,...)$. It is easy to check inductively that the resulting {\sc (km)} iterates are given by $x^n = (p^n_0,p^n_1,....,p^n_n,0,0,....)$,
where $p^n_k = \mathbb{P}(S_n= k)$ is the distribution of a sum $S_n=X_1+\cdots+X_n$ of independent Bernoullis  with $\PP(X_i=1)=\alpha_i$.

A well-known result by Darroch \cite{darroch1964} establishes that the distribution of $S_n$ is bell-shaped, from which it follows that
$$\|x^n-Tx^n\|_1 = 2\max_{0\le k\le n} p^n_k$$ 
the maximum being attained either at $k=\lfloor\mu\rfloor$ or $k=\lceil\mu\rceil$ (or both)
where $\mu=\alpha_1+\ldots+\alpha_n$.
Moreover,  taking $B_n(\bar\alpha)\sim \mbox{Binomial}(n,\bar\alpha)$ with $\bar\alpha = \frac{1}{n}\sum_{i=1}^n \alpha_i$,
a result from Hoeffding \cite{hoeffding1956distribution} shows that for $0\le b \le n\bar\alpha \le c\le n$ we have (see also Xu \& Balakrishnan \cite{xu2011convolution})
$$\mathbb{P}(b\le S_n \le c) \ge \mathbb{P}(b\le B_n(\bar\alpha) \le c).$$
Taking $b=\lfloor n\bar\alpha \rfloor$ and $c=\lceil n\bar\alpha \rceil$ it follows that 
\begin{align*}
2\max_{0\le k\le n} p^n_k & \ge p^n_{\lfloor n\bar\alpha \rfloor}+p^n_{\lceil n\bar\alpha \rceil} \\
& = \mathbb{P}(\lfloor n\bar\alpha \rfloor\le S_n \le \lceil n\bar\alpha \rceil) \\
& \ge \mathbb{P}(\lfloor n\bar\alpha \rfloor\le B_n(\bar\alpha) \le \lceil n\bar\alpha \rceil) .
\end{align*}

Let us define $f_n:[0,1]\rightarrow [0,1]$ by $f_n(x) = \mathbb{P}(\lfloor nx \rfloor\le B_n(x) \le \lceil nx \rceil)$. We want to compute the minimum value of $f_n$. 
Firstly, we observe that $f_n$ is symmetric with respect to $x=\frac{1}{2}$, i.e. $f_n(x) = f_n(1-x)$. 
Secondly, we note that $f_n$ is discontinuous at the points of the form $\frac{k}{n}$ for $k=1,...,n-1$. In fact, one can check that
 $f_n(\frac{k}{n})> f_n((\frac{k}{n})^-)\ge f_n((\frac{k}{n})^+)$ for all $k = 1,...,\lfloor \frac{n}{2}\rfloor$, and symmetrically $f_n(\frac{k}{n})> f_n((\frac{k}{n})^+)\ge f_n((\frac{k}{n})^-)$ for all $k = \lceil \frac{n}{2}\rceil,...,n-1$.  On each open interval $]\frac{k}{n},\frac{k+1}{n}[$ the function $f_n$ is differentiable and concave (see Figure \ref{fig:binomial})
and its infimum is attained asymptotically by approaching the extreme of the interval which is closest to $\frac{1}{2}$, namely 
$$x_{k,n}^*=\left\{\begin{array}{cl}
\frac{k+1}{n}&\mbox{if $\frac{k+1}{n}\leq \frac{1}{2}$},\\[0.5ex]
\frac{k}{n}&\mbox{if $\frac{k}{n}\geq \frac{1}{2}$},\\[0.5ex]
\mbox{both}&\mbox{if $\frac{k}{n}< \frac{1}{2}<\frac{k+1}{n}$.} 
\end{array}\right.
$$
 After some straightforward computations, one can conclude that  the infimum of $f_n$ over the full interval $[0,1]$
 occurs when $x$ tends to $\frac{1}{n}\lfloor \frac{n}{2}\rfloor$ from the right and/or when $x$ tends to 
 $\frac{1}{n}\lceil \frac{n}{2}\rceil$ from the left (see Figure \ref{fig:binomial}). 
 \begin{figure}[t]
\centering
\includegraphics[scale=0.5]{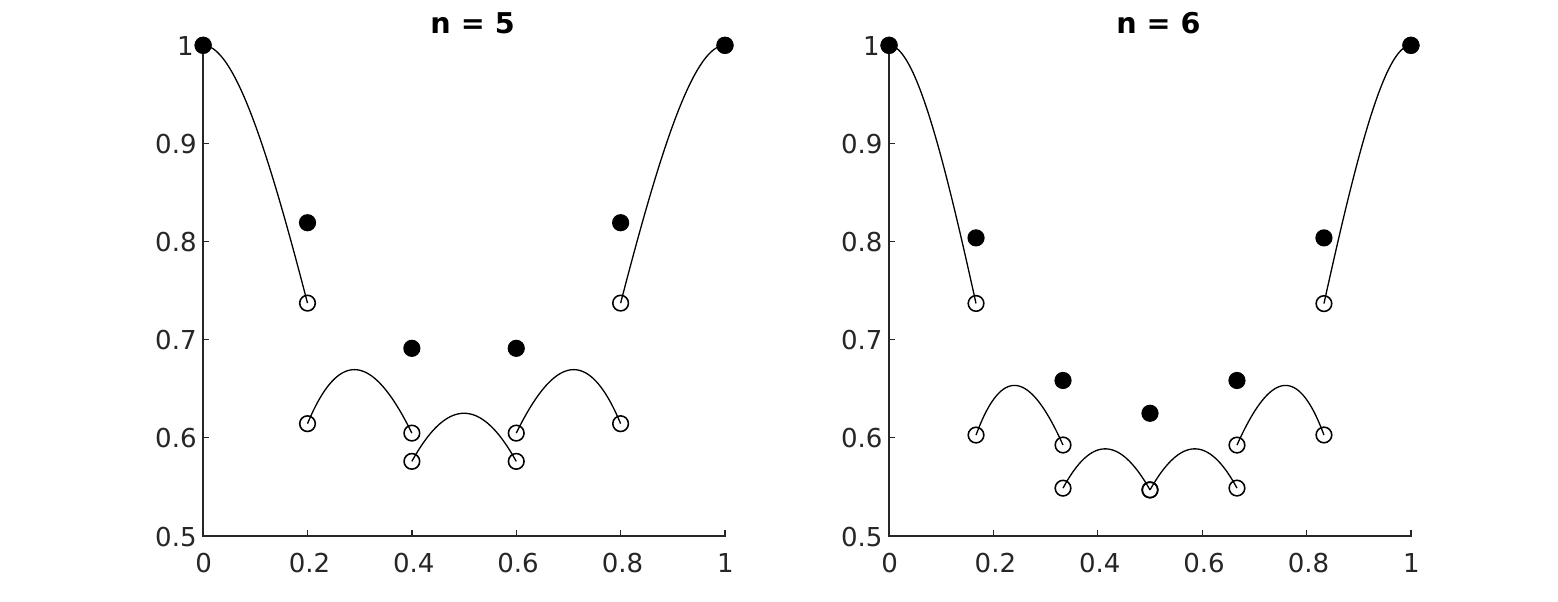}
\caption{The function $f_n$ for $n=5$ and  $n=6$.}
\label{fig:binomial}
\end{figure}

 In particular, when $n=2m$ is even, the infimum is obtained when 
 approaching $x=\frac{1}{2}$ either from the right or the left, with $\inf f_{2m} = \frac{2m+1}{m+1}\frac{1}{4^{m}} \binom{2m}{m}$. 
We observe that $\frac{1}{m+1}\binom{2m}{m}$ is a Catalan number, so that using the bound in Dutton \& Brigham \cite{dutton1986computationally} we  get
$$\mbox{$\|x^n-Tx^n\|_1 \ge \inf f_n\ge \frac{2m+1}{m+1}\sqrt{\frac{4m-1}{4m}}\frac{1}{\sqrt{\pi m}}\ge \frac{1}{\sqrt{n}}$.}$$ 
If $n = 2m+1$ is odd, then $\inf f_{2m+1}= \binom{2m+1}{m}\left(\frac{m(m+1)}{(2m+1)^2}\right)^m$. Using this expression along 
with the expression for the even case, it is easy to check that $\inf f_{2m+1}\ge \inf f_{2m+2}$, and therefore we conclude
$$\mbox{$\|x^n-Tx^n\|_1 \ge \inf f_{n} \ge \inf f_{n+1} \ge \frac{1}{\sqrt{n+1}}$}$$
completing the proof.
\end{proof}

\end{appendices}

{\bf Acknowledgements.} We thank professor Simeon Reich (Israel Institute of Technology) for his interest in this paper, his valuable comments, and for pointing out relevant references that helped us improve the introductory section. We also thank the two anonymous referees who carefully read the paper and provided important insights that helped us improve the presentation. The work of Juan Pablo Contreras was supported by a doctoral scholarship from ANID-PFCHA/Doctorado Nacional/2019-21190161. Roberto Cominetti
gratefully acknowledges the support provided by the research grant FONDECYT 1171501.


\begin{thebibliography}{37}
\providecommand{\natexlab}[1]{#1}
\providecommand{\url}[1]{\texttt{#1}}
\expandafter\ifx\csname urlstyle\endcsname\relax
  \providecommand{\doi}[1]{doi: #1}\else
  \providecommand{\doi}{doi: \begingroup \urlstyle{rm}\Url}\fi

\bibitem[Aronszajn and Panitchpakdi(1956)]{aronszajn1956extension}
 Aronszajn N. and  Panitchpakdi P.,
\newblock Extension of uniformly continuous transformations and hyperconvex
  metric spaces.
\newblock \emph{Pacific Journal of Mathematics}, 6(3):405--439, 1956.

\bibitem[Bailion et~al.(1978)Jean-Bernard Baillon, Bruck, and Reich]{bailion1978asymptotic}
Baillon J.B., Bruck R.E. and Reich S.,
\newblock On the asymptotic behavior of non-expansive mappings and semigroups in
  Banach spaces.
\newblock \emph{Houston Journal of Mathematics}, 4(1):1--10, 1978.

\bibitem[Baillon and Bruck(1992)]{baillon1992optimal}
Baillon J.B. and Bruck R.E.,
\newblock Optimal rates of asymptotic regularity for averaged non-expansive
  mappings.
\newblock \emph{World Scientific Publishing Co. Pte. Ltd., PO Box},
  128:\penalty0 27--66, 1992.

\bibitem[Baillon and Bruck(1996)]{baillon1996rate}
Baillon J.B. and Bruck R.E.,\newblock The rate of asymptotic regularity is $o(1/\sqrt{n})$.
\newblock In \emph{Theory and applications of nonlinear operators of accretive and
  monotone types}, \emph{Lecture Notes in Pure and Appl. Math.}, vol. 178:51--81. Dekker, New York, 1996.


\bibitem[Bauschke and Combettes(2011)]{bauschke2011convex}
Bauschke H.H. and Combettes P.L.,
\newblock \emph{Convex analysis and monotone operator theory in Hilbert
  spaces}, CMS Books in Mathematics, vol. 408.
\newblock Springer, 2011.

\bibitem[Berinde(2007)]{berinde}
Berinde V.,
\newblock \emph{Iterative Approximation of Fixed Points}, 
  {Lecture Notes in Mathematics}, vol. 1912.
\newblock Springer-Verlag Berlin Heidelberg, 2007.

\bibitem[Borwein et~al.(1992)Borwein, Reich, and
  Shafrir]{borwein1992krasnoselski}
 Borwein J.,  Reich S. and Shafrir I.,
\newblock Krasnoselski-Mann iterations in normed spaces.
\newblock \emph{Canadian Mathematical Bulletin}, 35(1):21--28, 1992.

\bibitem[Bravo and Cominetti(2018)]{bravo2018sharp}
 Bravo M. and Cominetti R.,
\newblock Sharp convergence rates for averaged non-expansive maps.
\newblock \emph{Israel Journal of Mathematics}, 227(1):163--188, 2018.

\bibitem[Bravo et~al.(2021)Bravo, Champion, and Cominetti]{bcc2020}
 Bravo M., Champion T. and Cominetti R.,
\newblock Universal bounds for fixed point iterations via optimal transport
  metrics.
\newblock \emph{arXiv:2108.00300v1}, pp. 1--21, 2021.

\bibitem[Browder(1967)]{browder1967convergence}
Browder F.E.,
\newblock Convergence of approximants to fixed points of non-expansive nonlinear
  mappings in Banach spaces.
\newblock \emph{Archive for Rational Mechanics and Analysis}, 24(1):82--90, 1967.

\bibitem[Browder and Petryshyn(1966)]{browder1966solution}
Browder  F.E. and Petryshyn W.V.,
\newblock The solution by iteration of nonlinear functional equations in Banach
  spaces.
\newblock \emph{Bulletin of the American Mathematical Society}, 72(3):571--575, 1966.

\bibitem[Colao and Marino(2021)]{colao2021}
Colao V. and Marino G.,
\newblock On the rate of convergence of Halpern iterations.
\newblock Preprint, pp. 1--8, 2021.

\bibitem[Cominetti et~al.(2014)Cominetti, Soto, and Vaisman]{cominetti2014rate}
Cominetti R., Soto J.A. and Vaisman J.,
\newblock On the rate of convergence of Krasnosel'ski{\u{\i}}-Mann iterations
  and their connection with sums of Bernoullis.
\newblock \emph{Israel Journal of Mathematics}, 199(2):757--772, 2014.

\bibitem[Darroch(1964)]{darroch1964}
Darroch J.N.,
\newblock On the distribution of the number of successes in independent trials.
\newblock \emph{The Annals of Mathematical Statistics}, 35(3):1317--1321, 1964.

\bibitem[Diakonikolas(2020)]{diakonikolas2020Halpern}
Diakonikolas J.,
\newblock Halpern iteration for near-optimal and parameter-free monotone
  inclusion and strong solutions to variational inequalities.
  \emph{Proceedings of Machine Learning Research}, 125:1--24, 2020.
  

\bibitem[Drori and Teboulle(2014)]{drori2014performance}
Drori  Y. and Teboulle M.,
\newblock Performance of first-order methods for smooth convex minimization: a
  novel approach.
\newblock \emph{Mathematical Programming}, 145(1-2):451--482, 2014.

\bibitem[Dutton and Brigham(1986)]{dutton1986computationally}
Dutton R.D. and Brigham R.C.,
\newblock Computationally efficient bounds for the Catalan numbers.
\newblock \emph{European Journal of Combinatorics}, 7(3):211--213, 1986.

\bibitem[Halpern(1967)]{Halpern1967fixed}
Halpern B.,
\newblock Fixed points of nonexpanding maps.
\newblock \emph{Bulletin of the American Mathematical Society}, 73(6):957--961, 1967.

\bibitem[Hoeffding(1956)]{hoeffding1956distribution}
Hoeffding W.,
\newblock On the distribution of the number of successes in independent trials.
\newblock \emph{The Annals of Mathematical Statistics}, 27(3):713--721, 1956.

\bibitem[Kim(2021)]{kim2021accelerated}
 Kim D.,
\newblock Accelerated proximal point method for maximally monotone operators.
\newblock \emph{Mathematical Programming}, 190(1):57--87, 2021.

\bibitem[Kohlenbach(2010)]{kohlenbach2010logical}
 Kohlenbach U.,
\newblock On the logical analysis of proofs based on nonseparable hilbert space
  theory.
\newblock In \emph{Proofs, Categories and Computations. Essays in Honor of Grigori
  Mints,} Eds. Solomon Fefferman and Wilfried Sieg, \emph{College Publications}, pp. 131--143, 2010.

\bibitem[Kohlenbach(2011)]{kohlenbach2011quantitative}
 Kohlenbach U.,
\newblock On quantitative versions of theorems due to F.E. Browder and R.
  Wittmann.
\newblock \emph{Advances in Mathematics}, 226(3):2764--2795, 2011.

\bibitem[K{\"o}rnlein(2015)]{kornlein2015quantitative}
K{\"o}rnlein D.,
\newblock Quantitative results for Halpern iterations of non-expansive mappings.
\newblock \emph{Journal of Mathematical Analysis and Applications},
  428(2):1161--1172, 2015.

\bibitem[Krasnosel'ski{\u{\i}}(1955)]{krasnosel1955two}
Krasnosel'ski{\u{\i}} M.A.,
\newblock Two remarks on the method of successive approximations.
\newblock \emph{Uspekhi Matematicheskikh Nauk}, 10:123--127,
  1955.

\bibitem[Leustean(2007)]{leustean2007rates}
 Leustean L.,
\newblock Rates of asymptotic regularity for Halpern iterations of non-expansive
  mappings.
\newblock \emph{Journal of Universal Computer Science}, 13(11):1680--1691, 2007.

\bibitem[Lieder(2021)]{lieder2021convergence}
 Lieder F.,
\newblock On the convergence rate of the Halpern-iteration.
\newblock \emph{Optimization Letters}, 15(2):405--418,
  2021.

\bibitem[Mann(1953)]{mann1953mean}
 Mann W.R.,
\newblock Mean value methods in iteration.
\newblock \emph{Proceedings of the American Mathematical Society}, 4(3):506--510, 1953.

\bibitem[Nesterov(2018)]{nesterov2018}
 Nesterov Y.,
\newblock \emph{Lectures on Convex Optimization},  \emph{Springer
  Optimization and Its Applications}, vol. 137.
\newblock Springer International Publishing, 2018.

\bibitem[Reich(1975)]{reich1975fixed}
 Reich S.,
\newblock Fixed point iterations of non-expansive mappings.
\newblock \emph{Pacific Journal of Mathematics}, 60(2):195--198, 1975.

\bibitem[Reich(1979)]{reich1979}
 Reich S.,
\newblock Weak convergence theorems for non-expansive mappings in Banach
  spaces.
\newblock \emph{J. Math. Anal. Appl.}, 67:274--276, 1979.

\bibitem[Reich(1980)]{reich1980}
 Reich S.,
\newblock Strong convergence theorems for resolvents of accretive operators in Banach spaces.
\newblock \emph{J. Math. Anal. Appl.}, 75:287--292, 1980.

\bibitem[Reich(1994)]{reich1994}
 Reich S.,
\newblock Approximating fixed points of nonexpansive mappings.
\newblock \emph{Panamerican Math. J.}, 4(2):23--28, 1994.

\bibitem[Ryu and Yin(2021)]{ryu2021large}
Ryu E.K. and  Yin W.,
\newblock Large-scale convex optimization via monotone operators. Book draft. 
To be published with \emph{Cambridge University Press}, 2022.

\bibitem[Sabach and Shtern(2017)]{sabach2017first}
Sabach S. and  Shtern S.,
\newblock A first order method for solving convex bilevel optimization
  problems.
\newblock \emph{SIAM Journal on Optimization}, 27(2):640--660, 2017.

\bibitem[Wittmann(1992)]{wittmann1992approximation}
 Wittmann R.,
\newblock Approximation of fixed points of non-expansive mappings.
\newblock \emph{Archiv der Mathematik}, 58(5):486--491,
  1992.

\bibitem[Xu(2002)]{xu2002iterative}
 Xu H-K.,
\newblock Iterative algorithms for nonlinear operators.
\newblock \emph{Journal of the London Mathematical Society}, 66(1):240--256, 2002.

\bibitem[Xu and Balakrishnan(2011)]{xu2011convolution}
Xu M.  and Balakrishnan N.,
\newblock On the convolution of heterogeneous Bernoulli random variables.
\newblock \emph{Journal of Applied Probability}, 48(3):877--884, 2011.

\end{thebibliography}
\end{document}